\documentclass{amsart}
\usepackage{times, mathrsfs, epsfig}
\usepackage{amsmath, amsthm, amssymb}
\input{diagrams}

\begin{document}

\newtheorem{thm}{Theorem}[section]
\newtheorem{lem}[thm]{Lemma}
\newtheorem{cor}[thm]{Corollary}
\newtheorem{conj}[thm]{Conjecture}
\newtheorem{qn}[thm]{Question}
\newtheorem{claim}[thm]{Claim}
\newtheorem{fact}[thm]{Fact}
\newtheorem{ex}[thm]{Exercise}

\theoremstyle{definition}
\newtheorem{defn}[thm]{Definition}
\newtheorem{note}[thm]{Notation}

\theoremstyle{remark}
\newtheorem*{rmk}{Remark}
\newtheorem{exa}[thm]{Example}

\def\square{\hfill${\vcenter{\vbox{\hrule height.4pt \hbox{\vrule width.4pt
height7pt \kern7pt \vrule width.4pt} \hrule height.4pt}}}$}

\def\R{\mathbb R}
\def\Z{\mathbb Z}
\def\CP{\mathbb {CP}}
\def\H{\mathbb H}
\def\h{\mathscr H}
\def\E{\mathbb E}
\def\S{\mathbb S}
\def\F{\mathscr F}
\def\D{\mathscr D}
\def\D3{{\Bbb D}^3}
\def\C{\mathbb C}
\def\Q{\mathbb Q}
\def\A{\mathbb A}
\def\til{\widetilde}
\def\M{\mathscr M}
\def\O{\mathscr O}
\def\L{\mathscr L}
\def\MH{\mathscr{MH}}
\def\MC{\text{MC}}
\def\N{\mathbb N}
\def\RP{\mathbb{RP}}
\def\T{\mathscr T}
\def\U{\mathscr U}
\def\V{\mathscr V}
\def\G{\mathscr G}
\def\gv{\mathfrak{gv}}
\def\P{\mathscr P}
\def\B{\mathscr B}
\def\X{\mathscr X}
\def\Aut{\text{Aut}}
\def\Inn{\text{Inn}}
\def\Out{\text{Out}}
\def\Isom{\text{Isom}}
\def\homeo{\text{Homeo}}
\def\stab{\text{stab}}
\def\hom{\text{Hom}}
\def\teich{\text{Teich}}
\def\u{{\text{univ}}}
\def\tr{\text{tr}}
\def\Li{\text{Li}}
\def\1{{\bf 1}}
\def\id{\text{id}}

\newenvironment{pf}{{\it Proof:}\quad}{\square \vskip 12pt}
\newenvironment{spf}{{\it Sketch of proof:}\quad}{\square \vskip 12pt}

\title{Problems in Foliations and Laminations of $3$--manifolds}

\author{Danny Calegari}
\address{Department of Mathematics \\ California Institute of Technology
\\ Pasadena CA 91125}
\email{dannyc@its.caltech.edu}
\date{9/8/2002. Version 0.78}
\maketitle

\section{Introduction and basic definitions}

\subsection{Notation}
I will try to use consistent notation throughout, so for instance, $\F,\G$ denote
foliations, $\Lambda$ denotes a lamination, $\lambda,\mu,\nu$ denote leaves, etc.
For a given object $X$ in a manifold $M$, $\til{X}$ will denote the pullback of $X$
to the universal cover $\til{M}$.

\subsection{Attribution}
I have tried to credit questions to their originators. There are
certain caveats to this, however:
\begin{itemize}
\item{I have made no attempt to verify that my sources were themselves
original.}
\item{I have occasionally edited or reformulated questions; in so doing,
it is possible that I have perverted the intent of the originator.}
\item{I have not systematically checked other problem lists to see if
the questions here were posed independently elsewhere.}
\end{itemize}

\subsection{Other problem lists}
The format and taste of this problem list is greatly influenced by
Gabai's well known problem list \cite{Gabprob}. However I have tried
to choose questions which complement, rather than
overlap with, Gabai's choice of problems.

\subsection{Basic objects}

A reference for the basic theory of foliations is \cite{CC_book}. An overview of
the state of the subject as it stood in 1990 is contained in \cite{Gab90}. \cite{Gabprob}
also contains background and numerous examples.

\begin{defn}
A codimension one foliation $\F$ of a $3$--manifold $M$ is {\em taut} if there is
a circle $\gamma$ transverse to $\F$ intersecting every leaf.
\end{defn}

A basic reference is \cite{Gab83}. Note that a foliation of an atoroidal $3$--manifold
is taut iff it contains no torus leaf.

For a taut foliation $\F$ of $M$, $\til{M}$ may be taken to be an open subset of $\R^3$
so that $\til{\F}$ consists of the intersection of $\til{M}$ with the horizontal planes in
$\R^3$. These leaves are all disks, and the leaf space $L$
of $\til{\F}$ is a simply--connected $1$--manifold. There is a canonical {\em holonomy
representation}
$$\rho_H:\pi_1(M) \to \homeo(L)$$

\begin{defn}
A taut foliation is said to be or have
\begin{itemize}
\item{{\em $\R$--covered} if the leaf space of $\til{\F}$ is homeomorphic to $\R$.}
\item{{\em one--sided branching} if the leaf space of $\til{\F}$ branches in one direction.}
\item{{\em two--sided branching} if the leaf space of $\til{\F}$ branches in both directions.}
\end{itemize}
\end{defn}

All possibilities are realized, for foliations of atoroidal
$3$--manifolds.

\begin{defn}
A codimension one {\em lamination} $\Lambda$ in a $3$--manifold $M$ is a
foliation by two--dimensional leaves of a closed subset of $M$. A lamination is {\em essential}
if it contains no sphere leaf or torus leaf bounding a solid torus, if complementary
regions are irreducible and have incompressible boundary, and admit no compressing
monogons.
\end{defn}

A taut foliation is an example of an essential lamination. A closed union of
leaves of a taut foliation is another example.

\begin{defn}
The complementary regions to an essential lamination fall into two
kinds of pieces, which can be chosen uniquely up to isotopy:
$I$--bundles $\Sigma_i \times I$ over noncompact surfaces $\Sigma_i$,
called {\em interstitial regions}, and compact pieces
called {\em guts} $\mathfrak{G}_i$ meeting
the $I$--bundles along annuli called {\em interstitial annuli}. An essential lamination
is {\em genuine} if it has nonempty gut --- i.e. if some complementary region is not
an $I$--bundle.
\end{defn}

To construct the partition of a complementary region $C$ into guts and
interstices, take as a first approximation of the interstices the
characteristic $I$--bundle of $C,\partial C$, and then remove all
$I$--bundles over non--compact components.

For more details on essential and genuine laminations, see \cite{GO89}, \cite{GK93} and \cite{GK98}.
Note that though foliations can always be co--oriented in some finite cover of the ambient
manifold, laminations admit local obstructions to co--orientability.

\begin{defn}
A genuine lamination has {\em solid torus guts} if all gut regions are
neutered ideal polygon bundles over $S^1$. It is
{\em full} if some complementary region is an ideal polygon
bundle over $S^1$. It is {\em very full} if {\em every} complementary region is a bundle region.
\end{defn}

There is a procedure to turn {\em co--oriented} very full genuine laminations into
taut foliations, by filling in the complementary regions with {\em monkey saddles}. These
are disks which are asymptotic along the boundary leaves of a complementary region
in alternating directions. The complementary regions can be given a bundle structure
over $S^1$ where these saddle leaves are the fibers.

\begin{defn}
A {\em pseudo--Anosov flow} $X$ on a manifold $M$ is a flow which, away from a finite
number of closed orbits $\gamma_i$, is {\em Anosov} --- that is, there is a decomposition
of $TM$ as a sum $E_s \oplus E_u \oplus TX$ where the time $t$
flow contracts $E_u$ by some
factor $\lambda^{-t}$ and expands $E_s$ by $\lambda^t$, for some $\lambda>0$. Along the
singular orbits, the local picture is that of a {\em semi--branched cover}
of an Anosov orbit of degree $n/2$ for $n \ge 3$. The
distributions $E_s$ and $E_u$ are integrable, and integrate to give two singular
foliations $\F_s,\F_u$ the {\em stable} and {\em unstable} integral foliations. If
one splits open $\F_s,\F_u$ along the singular leaves and takes the closure, one obtains
very full genuine laminations $\Lambda_s,\Lambda_u$.
\end{defn}

The basic examples of pseudo--Anosov flows are the suspensions of pseudo--Anosov
automorphisms of surfaces. See \cite{Thu88}. 

\begin{defn}
A lamination $\Lambda$ has {\em hyperbolic leaves} if the leaves of $\til{\Lambda}$
are (uniformly) quasi--isometric to $\H^2$ with the path metric induced from some
metric on $M$.
\end{defn}

A theorem of Candel (\cite{Can93}) implies that for $M$ atoroidal, any essential
lamination $\Lambda$ can be given a metric of {\em constant} curvature $-1$. In
particular, all such laminations have hyperbolic leaves.

\begin{defn}
Suppose $\F$ is a taut foliation of $M$ with hyperbolic leaves. A {\em universal circle}
for $\F$ is a circle $S^1_\u$ with a group action
$$\rho_\u:\pi_1(M) \to \homeo(S^1_\u)$$ satisfying the following
properties:
\begin{enumerate}
\item{For all leaves $\lambda$ of $\til{\F}$ there is a {\em monotone} map
$$\phi_\lambda:S^1_\u \to S^1_\infty(\lambda)$$
where $\phi_\lambda$ varies continuously as a function of $\lambda$.
Here a map between circles is {\em monotone} if the preimage of every
point is contractible.}
\item{For each leaf $\lambda$ of $\til{\F}$ and each 
$\alpha \in \pi_1(M)$ the following diagram commutes:
\begin{diagram}
S^1_\u                          & \rTo^{\rho_\u(\alpha)}         & S^1_\u \\
\dTo^{\phi_{\lambda}} &                         & \dTo^{\phi_{\alpha(\lambda)}} \\
S^1_\infty(\lambda)   & \rTo^{\alpha}         & S^1_\infty(\alpha(\lambda)) \\
\end{diagram}}
\item{For any leaf $\lambda$, the {\em associated gaps} are the
maximal connected open intervals in $S^1_\u$ mapped to points by
$\phi_\lambda$. The complement of the gaps in $S^1_\u$ is the
{\em core} associated to $\lambda$. Then for any pair $\mu,\lambda$
of incomparable leaves, the core associated to $\lambda$ is
contained in a single gap associated to $\mu$, and vice versa.}
\end{enumerate}
\end{defn}

Universal circles are introduced in \cite{Thu??}; examples and discussion are
found in \cite{Cal00}, \cite{Calp1} and \cite{CalPhD}. Since the map between the
universal circle and $S^1_\infty(\lambda)$ is monotone for each $\lambda$, it easily
follows that the Euler class of the representation $\rho_\u$ is equal to the Euler
class of $T\F$.

\begin{defn}
A {\em branched surface} $\B$ is a certain kind of $2$--complex in a $3$--manifold with a
well--defined $C^1$ tangent space everywhere. The place where the branched surface is
not a surface is called the {\em branch locus}; it is a $1$--complex, and the
branch locus is not a $1$--manifold only at isolated points, where two lines of
the branch locus cross each other transversely.
\end{defn}

Along the nonsingular part of the branch locus of a branched surface,
there is a natural transverse orientation,
which points in the direction in which the two branches come together. Where two
arcs of the branch locus cross, there is one complementary quadrant where both
co--orientations are inward, and one where both co--orientations are outward.

\begin{defn}
A {\em sink disk} for a branched surface is an embedded disk $D \subset \B$
with boundary contained in the branch locus, such that the transverse orientation
points inward all along $\partial D$.
\end{defn}

Sink disks are introduced by Li in \cite{Li00}. He shows, amongst other things, that
branched surfaces without sink disks fully carry essential laminations, and that
every essential lamination (with some very simple exceptions) is carried by a branched
surface without sink disks.

\subsection{$3$--manifold topology}

Throughout this paper, numerous reference is made, implicitly or explicitly, to
elements of the theory of $3$--manifolds. A basic reference for this theory
is \cite{Hem76}. A reference for the theory of geometrization is \cite{Thu_notes}.

\subsection{Acknowledgements}

Many thanks to Charles Delman, Tao Li, Curt McMullen, Saul Schleimer and participants
in the problem session on foliations and laminations at the 2001 Georgia International
Topology Conference. Ian Agol and Nathan Dunfield should be singled out for
special thanks, for their considerable contributions to this problem list. 

\section{Existence questions}

\begin{qn}
Which hyperbolic $3$--manifolds admit taut foliations? Give an effective
procedure to decide if a hyperbolic $3$--manifold admits a taut foliation.
For a useful property $\sf{p}$ of a $3$--manifold, give an explicit
construction of infinitely many $3$--manifold with
property $\sf{p}$ with/without a taut foliation. Same question for essential
laminations.
\end{qn}
\begin{rmk}
\begin{enumerate}
\item{It is known by \cite{RSS} that there exist infinitely many
hyperbolic $3$--manifolds
without taut foliations. These manifolds are so special that their fundamental
groups admit no nontrivial actions on (potentially non--Hausdorff) $1$--manifolds. Here
an action on a non--Hausdorff $1$--manifold is trivial if it has a global
fixed point {\em on the maximal Hausdorff quotient}.
The proof is somewhat {\it ad hoc} and relies on the special
structure of the fundamental groups of these manifolds (they are Dehn
fillings of certain punctured torus bundles with negative trace). It should be
remarked that Hatcher initially speculated that such manifolds would be
good candidates to be taut foliation free.

The method of proof is not entirely satisfactory for two reasons.
Firstly, the {\em certificate} that the manifolds in question admit no taut
foliation is very long, and it seems unlikely that it could be further compressed within a
similar approach. (See \cite{CD01} for a different method which provides a
much shorter certificate that the Weeks manifold admits no taut foliation). Secondly,
the criterion by which the existence of taut foliations is ruled out is
not {\it a priori} sharp. For instance, if $f:M \to N$ has positive
degree and $N$ admits a taut foliation $\F$, then $\pi_1(M)$ acts on
the leaf space of $\til{\F}$ without a global fixed point; on the other hand,
it is unclear whether such an $M$ admits a taut foliation in general. In particular, it
is important to understand when a $3$--manifold admitting a nontrivial action on a
non--Hausdorff $1$--manifold admits a taut foliation.

The methods of \cite{RSS} are
extended by Fenley in \cite{Fen_nolam} to show that infinitely many of these
hyperbolic manifolds actually admit no essential laminations, by showing that
their fundamental groups admit no nontrivial actions on order trees. It suffers from
the same shortcomings as \cite{RSS}, and is therefore perhaps more of theoretical
and philosophical than of practical interest.}
\item{Work of Agol \cite{Agol} gives strong restrictions on the topology of the guts
of a tight genuine lamination on a hyperbolic manifold $M$ of small volume.
If a taut foliation $\F$ is not monotone equivalent to a minimal foliation,
some minimal set is a genuine lamination. Thus it might be interesting
to give examples of $3$--manifolds without taut foliations of large volume
and/or injectivity radius.}
\item{Agol and Li \cite{AL02}, combining work of Li \cite{Li00} and Brittenham \cite{Bri95},
give an algorithm to decide whether a given $3$--manifold admits an essential
lamination. For, if a $3$--manifold admits a nowhere dense essential lamination,
Brittenham shows it admits an essential
lamination in normal form with respect to any triangulation, and therefore is
carried by one of a finite constructible set of branched surfaces. With more
work, using technology of Gabai \cite{Gab99}, one can show that every
essential nowhere dense lamination is {\em fully carried} by one of a finite
constructible set of branched surfaces. Li shows that
every essential lamination can, after blowing up, be carried by a branched surface
without a sink disk (except for the special case of a foliation of $T^3$ by planes).
So the algorithm proceeds by taking each of the finite constructed set of branched
surfaces, and splits it open repeatedly in all possible ways. After some finite
time, either every branched surface reveals an obstruction to being fully split
open, or some branched surface can be split open to a branched surface with no
sink disk. How computationally efficient is this algorithm? Is there any case
in which it might actually be implemented on a computer, e.g. for the Seifert--Weber
dodecahedral space? 

Note that the method of Agol and Li decides
whether a given $3$--manifold admits a taut foliation. For, one can check that
one of the finite constructible set of branched surfaces which fully carries an
essential lamination has $I$--bundle complementary region.}
\item{Properties $\sf{p}$ for which this question is interesting might include
Hakenness, or more generally, laminar properites. Obstructions to
the existence of a taut foliation which are homological in nature probably
do not obstruct the existence of genuine laminations. In particular, the
local non co--orientability is a strong advantage. Since taut foliations
can be perturbed to symplectically semi--fillable contact structures by
Eliashberg and Thurston \cite{ElTh}, one might be able to use contact
or symplectic $4$--manifold techniques (e.g. gauge theory) to rule out
the existence of taut foliations on certain manifolds which contain genuine
laminations. Conversely it would be
very interesting to give an example of a hyperbolic $3$--manifold with a
semi--fillable contact structure but no taut foliation.}
\item{For a given manifold $M$ and a hyperbolic knot $K \subset M$, the
set of surgeries on $K$ giving rise to a tautly foliated or laminar manifold, if
not everything, is probably quite complicated in some cases. There are
probably arithmetic criteria (e.g. some persistent lamination might be
co--orientable or not after certain surgeries and therefore fillable or
otherwise), inequalities (e.g. branched surfaces supporting
foliations realizing certain slopes), bounded cohomology (e.g. Milnor--Wood),
etc. There are existence results of Gabai--Mosher \cite{Mos} which say that
if $K$ is a hyperbolic knot,
there is a degeneracy slope $\delta$ associated to $K$ such that the stable and
unstable laminations coming from a pseudo--Anosov flow on $M-K$ stay essential
for surgery on all slopes $\alpha$ with $\Delta(\delta,\alpha)>1$. Bootstrap
techniques to determine degeneracy slopes are exploited in \cite{CD01} to
give a partial table of which small volume hyperbolic manifolds are laminar.}
\end{enumerate}
\end{rmk}

\begin{qn}
Is there an effective algorithmic
procedure to produce and recognize a hyperbolic
knot of depth $n$ for any given $n$? What about $\ge n$?
\end{qn}

\begin{rmk}
\begin{enumerate}
\item{Satellite knots of arbitrary depth were produced by Cantwell and Conlon.
A single explicit example of a hyperbolic knot of
depth $\ge 2$ is constructed by Kobayashi in \cite{Kob??}.
Brittenham has a class of examples which are conjecturally
of depth $\ge 2$. But what about higher depth? See also
question~\ref{dimension_depth}.}
\item{Algorithms exist to detect whether a knot is of depth $0$ (\cite{Hem})
or of depth $1$ but they are not very useful in practice.
Part of the problem is that finding Thurston norm minimizing
surfaces is NP complete by \cite{AHT}. On the other
hand, perhaps there is a fast algorithm to {\em check}
that a foliation is of minimal depth. A ``good'' algorithm
producing a knot of depth $n$ would also produce
a certificate for this fact. The paper \cite{Shi97} may be relevant here.}
\item{Lackenby's technique of taut ideal triangulations \cite{Lacktit} is useful in
practice for finding and especially {\em certifying} minimal genus representatives
for relative homology classes in open manifolds with torus boundary.}
\end{enumerate}
\end{rmk}

\begin{qn}
Given a collection $\mathscr{C}$
of topological or geometric types of surface, what $3$--manifolds
admit a taut foliation $\F$ whose leaves are all homeomorphic or coarsely
quasi--isometric to an element of the collection $\mathscr{C}$?
\end{qn}
\begin{rmk}
\begin{enumerate}
\item{Every topological surface is a leaf of a foliation (even a taut
foliation) of some $3$--manifold $M$, by a construction of Cantwell
and Conlon \cite{CaCo87}.}
\item{The Reeb stability theorem implies that a foliation of a $3$--manifold
with a spherical leaf is a bundle (possibly over an orbifold). Moreover,
taut foliations of hyperbolic $3$--manifolds do not contain torus or Klein bottle
leaves.}
\item{By a theorem of Candel \cite{Can93}, if $\F$ is a taut foliation of an
atoroidal $3$--manifold $M$, there is a metric on $M$ with respect to which every
leaf has constant curvature $-1$. In particular, the universal covers of all the leaves are
isometric to $\H^2$. Moreover, if there is a leaf (of $\F$, not $\til{\F}$)
with amenable growth, then its closure supports an invariant transverse measure,
and therefore contains leaves of polynomial growth, by Plante's theorem. In particular,
if $M$ is a rational homology sphere, every leaf of an orientable/co--orientable
taut foliation $\F$ has {\em uniformly} exponential growth.}
\item{In great generality this question is probably not so interesting, but
in some special cases it is probably quite interesting. For instance,
if $\mathscr{C}$ consists just of a plane, $M$ must be $T^3$ and the foliation
must be monotone equivalent to one by planes of irrational slope. If $\mathscr{C}$ consists
of just planes and cylinders and $\F$ is $\R$--covered, either $M$ is
Solv, or there are cylinders whose nontrivial loops are homotopic to each
other. If one adds the geometric constraint that the shortest loops in the
cylinders are of bounded length, if the ambient manifold is atoroidal and
therefore has word--hyperbolic fundamental group, these loops must be
homotopic by an annulus of controlled diameter etc.}
\item{The character of this problem would be very different if one asked for
foliations with the property that {\em almost all} leaves were homeomorphic or
coarsely quasi--isometric to an element of $\mathscr{C}$. If $\mu$ is a harmonic measure
on a foliated space, then any leaf $\lambda$ in the support of $\mu$ is quasi--homogeneous;
that is, any local geometric feature of $\lambda$ on any scale
recurs with definite density in
any sufficiently large subset of $\lambda$. This is more or less a consequence of
the fact that the support of $\mu$ is roughly the subset where leafwise random
walks are recurrent. In this vein, there is a very satisfying theorem of Ghys
that for $\Lambda$ {\em any} (abstract) compact Riemann surface lamination, and $\mu$
any harmonic measure, then $\mu$--almost all leaves are homeomorphic to one of the
following six possibilities:
\begin{itemize}
\item{The plane}
\item{The Loch--Ness monster}
\item{The cylinder}
\item{Jacob's ladder}
\item{$S^2 - \text{ Cantor set }$}
\item{$S^2 - \text{ Cantor set }$ with countably many handles accumulating to every point
in the removed set.}
\end{itemize}
Ghys' theorem is proved in \cite{eG95}.}
\end{enumerate}
\end{rmk}

\begin{qn}
Let $X$ be a vector field on a $3$--manifold. When is there a foliation $\F$ of $M$
transverse to $X$?
\end{qn}
\begin{rmk}
\begin{enumerate}
\item{If $M$ is Seifert--fibered and $X$ is a vector field tangent to the circle fibers,
there are many restrictions on the existence of $\F$. The {\em Milnor--Wood} inequality
says that the Euler class of a foliated circle bundle is at most as large in absolute
value as the Euler characteristic of the base orbifold \cite{Mil58},\cite{Woo69}. This
result can be phrased in terms of bounded cohomology, as saying that the norm of
the Euler class is equal to $1$ \cite{Gro82},\cite{Thu97}. Moreover,
results of Brittenham and Thurston \cite{Bri93}\cite{Thu72PhD} imply that any
taut foliation of a Seifert--fibered manifold is either transverse to the circle fibers
or contains ``vertical'' leaves --- leaves foliated by circles. Complete necessary and
sufficient conditions for the existence of a foliation transverse to the foliation
by circles of a Seifert fibered $3$--manifold were established by Naimi \cite{Nai}
building on results of Jankins and Neumann \cite{JN85}.}
\item{By analogy with the circle bundle case, one can ask whether
there exists a {\em slithering}
of $M$ over $S^1$ for which the slithering map $Z$ is realized by time $1$ flow along $X$.
Define an equivalence relation $\sim$ on $M$ which identifies pairs of points which
are an integer distance apart on a flow line of $X$ (parameterized by the flow speed). Then
one can think of $M/\sim$ as a circle bundle over the non--commutative space $M/X$, the leaf
space of $X$. One would like to prove an analogue of the Milnor--Wood
inequality for some appropriate fundamental
classes for $M/\sim$ and $M/X$. See \cite{Con94}. One can vary the question by allowing
(nonconstant) reparameterizations of $X$.}
\item{More generally, one can study a homology theory closely related to
the {\em Godbillon homology} of $X$, whose $i$--chains are equivalence classes of singular
maps of simplices to $M$, where two such maps $\sigma_1,\sigma_2:\Delta \to M$
are equivalent if there is a homotopy $S:\Delta \times I \to M$
where $S(\cdot,0) = \sigma_1(\cdot)$,
$S(\cdot,1) = \sigma_2(\cdot)$ and $S(p,I)$ is contained in a flowline of $X$. What is the
cohomology of the {\em bounded} cochains on this complex?}
\end{enumerate}
\end{rmk}

\section{Rigidity and Moduli}

\begin{qn}\label{polyhedral_pseudos}
Let $M$ be atoroidal.
Is there a natural refinement of the polyhedral structure of the unit ball of the
Thurston norm to a polyhedron $\P_\F$ whose faces parameterize the set of taut of
foliations of $M$ in the following sense:
\begin{enumerate}
\item{Each open face $c$ of $\P_\F$ corresponds to a
pseudo--Anosov flow $X_c$ or equivalence classes of pseudo--Anosov flows.
The taut foliations ``parameterized'' by $\P_\F$
can be isotoped to be transverse
(or {\em almost transverse}) to $X_c$, or some
sublamination is monotone equivalent
into the stable or unstable singular foliation of $X_c$.}
\item{Every geometric limit of taut foliations $\F_i$ associated to a cell $c$
should be associated to some cell of the closure $\overline{c}$.}
\item{There should be a natural polyhedral map from $\P_\F$ to (some polyhedral
subdivision of) the unit ball of the Thurston norm.}
\end{enumerate}
\end{qn}

\begin{rmk}
\begin{enumerate}
\item{Basic references for the Thurston norm are \cite{Thu86} and \cite{Oer86}.
For depth 0 (surface bundles), fibered homology classes are a union of top
dimensional faces, by \cite{Thu86}.
For depth 1, a refinement of the polyhedral structure
is necessary, but is accomplished by Cantwell and Conlon in (\cite{CCgeo}).
Their procedure applies to depth 1 foliations of
sutured manifolds, and therefore
inductively, should give a polyhedral parameterization of
all finite depth foliations
on a fixed manifold, up to some appropriate equivalence relation.}

\item{For a given fixed manifold $M$, Gabai has shown in \cite{Gab99}
there is a triangulation $\tau$ such that
every taut foliation of $M$ can be put into normal form with respect to $\tau$.
There is a natural polyhedral structure $\P_\tau$
on the space of all normal surfaces carried
by $\tau$. $\P_\F$ should be
intermediate between the unit ball for the Thurston norm and $\P_\tau$.

Note for certain classes of foliations
($\R$--covered, one--sided branching, minimal
foliations etc.) there are {\em existence}
results for transverse pseudo--Anosov flows
and, in many cases, stability results which
say that these flows are unique up to
isotopy for small deformations of the foliations (\cite{Cal00},\cite{Calp1},\cite{Calp3}) .
It seems plausible that a deformation
of normal foliations which preserves the normal
disk types in each simplex should
be mild enough to ensure stability of a transverse flow.
Note that for {\em finite depth foliations} one has the theorem of Lackenby 
that a taut sutured manifold hierarchy can be placed in a position
resembling normal form, relative to {\em any} triangulation (\cite{Lackexc}).}

\item{Kronheimer and Mrowka (\cite{KM97}) found a
characterization of the unit ball
of the dual Thurston norm by means of differential geometry. Is there some way
to detect the depth of a taut foliation or the
codimension of the face it corresponds to
by geometric properties of an approximating tight contact structure?}

\item{By Gabai and Mosher, there is a pseudo--Anosov flow almost transverse
to any finite depth foliation of a hyperbolic $3$--manifold \cite{Mos}.
For depth $>1$, there are
choices involved in the construction of the flow; any natural equivalence
relation on pseudo--Anosov flows would have to take this ambiguity into
account. Mosher in particular has developed a sophisticated picture of the
relationship between pseudo--Anosov and, more generally, $pA$ flows on $3$--manifolds
and the Thurston norm, for instance in the papers \cite{Mos91},\cite{Mos92a},\cite{Mos92b}.}

\item{Thurston's theorem that euler classes of taut foliations are
contained in the unit ball of the dual Thurston norm follows from the
existence of a {\em universal circle} for a taut foliation, and from
the Milnor--Wood inequality. One potential source of polyhedra $\P_\F$
might be unit balls of $L_1$ norms on certain natural quotients of
the vector space $\Omega^1(M)$ of $1$--forms on $M$.
Even an infinite dimensional
polyhedron parameterizing taut foliations would be very interesting.}
\end{enumerate}
\end{rmk}

\begin{qn}[McMullen]\label{Teichmuller_polynomial}
Generalize the {\em Teichm\"uller polynomial} from the fibered faces of the
Thurston norm ball to the other faces (of some possibly generalized polyhedron,
perhaps the polyhedron sought in question~\ref{polyhedral_pseudos}).
\end{qn}
\begin{rmk}
\begin{enumerate}
\item{A reference for the Teichm\"uller polynomial associated to a fibered
face of the Thurston norm ball is \cite{McMfib}. The fibered
face itself is the Newton polygon of the polynomial associated to that face;
amongst other information, the polynomial encodes the
(topological) entropy of the (topologically unique)
transverse pseudo--Anosov flow associated to the face. If the singular stable and
unstable foliations of the pseudo--Anosov flow are co--orientable, similar information
is encoded in Fried's dynamical zeta function \cite{Fr82}.}
\item{The entropy of the pseudo--Anosov flow arising from a surface bundle over a
circle is intimately related to the entropy of the automorphism of the fiber. Loosely
speaking, the Teichm\"uller polynomial is a kind of ``characteristic polynomial'' for
the action of the monodromy on the space of invariant transverse measures for an
invariant train track. \cite{Mat87a} gives a nice exposition of these ideas.}
\item{For a depth $n$ foliation where $n>0$, one has invariant laminations for the
end--periodic automorphisms of the fibers of greatest depth (see e.g. \cite{FenPhD}),
and one can study the entropy of these automorphisms. Global genuine laminations
which can be collapsed to pseudo--Anosov flows
are built up from these laminations at lowest depth inductively.
A transverse pseudo--Anosov flow $X$ for a fibration or a transversely
measured foliation $\F$ is like a Teichm\"uller geodesic for the family
of hyperbolic metrics on the leaves of $\F$ (see \cite{McMfib}).
The action of $\pi_1(M)$ on this Teichm\"uller geodesic
is conjugate to the action on the transverse measures in the appropriate projective class
for the singular foliations associated to $X$. For more general (not transversely measured)
foliations with transverse pseudo--Anosov flows, there might exist refined
polynomials with values in $\Z[GL(n,\R)]$
for some appropriate representation of $GL(n,\R)$ in $\homeo(\R)$. For instance, there
are many interesting slitherings of hyperbolic $3$--manifolds
with transverse $\til{SL(2,\R)}$ structure, associated
to branched covers over unit tangent bundles of hypebolic
orbifolds, where the branch loci are timelike with respect to the usual Lorentz structure
on $PSL(2,\R)$. See \cite{Thu97} or questions~\ref{SL2R_laminations}
and~\ref{automatic_order}.}
\end{enumerate}
\end{rmk}

\section{Minimal surfaces}

\begin{qn}\label{minimal_surface}
Suppose $\F$ is a taut foliation of $M$. Characterize the space of metrics on $M$ for which
$\F$ can be isotoped to consist of minimal surfaces.
\end{qn}
\begin{rmk}
\begin{enumerate}
\item{For general results on minimal surfaces relevant to this question, see \cite{Has86},
\cite{Pit81}, \cite{Rub01}.}
\item{By Sullivan \cite{Sul79} a $C^2$ foliation of a $3$--manifold can be made
minimal with respect to {\em some} Riemannian metric on the ambient manifold iff it is
taut. This is equivalent to the existence of a nonsingular divergence free {\em normal
vector field} --- i.e. one which is perpendicular to $T\F$ everywhere. Given any
transverse nonsingular vector field $X$ for which $\L_X \omega = 0$ for some nonsingular
$3$--form $\omega$, one can find a metric on $M$ for which $X$ is normal and $\omega$ is
the volume form.}
\item{A minimal surface in a hyperbolic manifold has intrinsic curvature $\le -1$ everywhere.
By Gauss--Bonnet, this puts an upper bound on the area of a compact minimal surface. Suppose
$\F$ is a taut foliation in a hyperbolic $3$--manifold enjoying some comparable bound
on the leafwise intrinsic curvature. If $\F$ admits a nontrivial invariant transverse measure
$\mu$ (not necessarily of full support), one obtains nontrivial constraints on
the growth rate and topology of leaves in the support of $\mu$.}
\item{Ben Andrews has constructed an area--minimizing
flow on surfaces with principle curvatures $<1$ to a minimal surface with the same properties.
Thurston has observed that elliptic PDE flows generally fail to preserve integrability of
plane fields, since after a short time, such flows tend to make structures real analytic;
but there are many $C^\infty$ foliations, even taut foliations, which admit no real
analytic structure, and even manifolds which admit no (nontrivial) transversely 
real analytic foliation at all. Even so, it might be interesting to understand the behaviour of
the leaves of $\F$ under Andrews' flow.}
\item{A minimal surface in $\E^3$ or $\H^3$ has a natural real analytic structure, which
one can think of as a kind of {\em tangential rigidity} condition. Zeghib shows in
\cite{aZ00} that such tangential rigidity for foliations implies a transverse
Lipschitz regularity; that is, any minimal foliation of $\E^3$ or $\H^3$ should be transversely
Lipschitz. This already shows that certain taut foliations of hyperbolic $3$--manifolds
cannot be made minimal with respect to the hyperbolic metric, and perhaps can be used to
show no foliation of a hyperbolic $3$--manifold by minimal surfaces can occur. Notice
in codimension $2$ a minimal foliation would actually have to be geodesic; Zeghib showed
that there are no nontrivial geodesic foliations of any codimension of any closed
hyperbolic manifold (\cite{aZ93}). On the other hand, if $\F$ is proper and all leaves
have the continuous extension property, one can inductively minimize leaves of fixed depth.
The closure of the union will be an essential (but not genuine) lamination, monotone equivalent
to the original foliation; but there seems to be no reason to expect it to be a foliation.}
\item{One weakened version of the question is as follows: is there an {\it a priori}
upper bound on the infimum of the mean curvature of a taut foliation $\F$ in a hyperbolic
$3$--manifold $M$, independent of $\F$ and $M$?}
\end{enumerate}
\end{rmk}

\begin{qn}
Given a collection of taut foliations $\F_i$ of $M$, what are the obstructions to finding
a metric on $M$ for which the $\F_i$ (after an isotopy) are simultaneously
minimal? 
\end{qn}
\begin{rmk} 
\begin{enumerate}
\item{Some kind of ``stably minimal'' property holds at the PL level.
For foliations $\F_i$ whose tangent fields are sufficiently close as $2$--plane fields,
we can find a triangulation $\tau$ of $M$ with respect to which the $\F_i$ can all
be made simultaneously normal, in such a way that the co--orientations on the
one skeleton of $\tau$ coming from the co--orientations on $\F_i$ agree.
For sufficiently close $\F_i$, we can find a triangulation such that the
$1$--skeleton with this orientation is recurrent; such a $1$--skeleton 
admits a combinatorial weight $w: \tau^1 \to \Z^+$ such that for each
vertex, the sum of the weights of the incoming edges is equal to the sum of
the weights of the outgoing edges. With respect to such a weight, these
foliations are ``minimal'' polyhedral surfaces. For similar ideas,
see \cite{Has86} and \cite{Gab??}.}
\item{Maybe this question become easier if we can replace the $\F_i$ by
monotone equivalent laminations $\Lambda_i$. This is especially plausible
if the $\F_i$ are finite depth --- the leaves at each depth could be
constructed inductively, and one could ignore the problem of filling up
the remaining gaps with the leaves of greatest depth.}
\end{enumerate}
\end{rmk}

\section{Reeb components}

Every closed $3$--manifold admits a foliation
with at most one Reeb component. One way to see
this is to use an open--book decomposition. A more
profound construction due to Gabai
shows the core of the Reeb component can be essentially
any knot (\cite{Gab83}). Every open $3$--manifold admits a foliation
with possibly infinitely many Reeb components. This
follows from Thurston's local construction of a
foliation from an anti--orientation on a triangulation (\cite{Thu73}),
and Moise's theorem that $3$--manifolds can be triangulated.

\begin{qn}
How many Reeb components must a foliation of an open $3$--manifold contain?
\end{qn}

\begin{rmk}
\begin{enumerate}
\item{There are two questions here, a topological
one and a geometric one. The topological
question is perhaps not so interesting. To make the geometric question precise,
some kind of quantification is necessary.
First, the $3$--manifold should have uniformly bounded geometry,
and then one should uniformly bound the local geometry of the foliation, and
either require all Reeb components to have bounded diameter, or weight them by
the length of a core circle.

Here it would be useful to distinguish between answers {\em none,
finitely many, infinitely many}, and if the third
alternative, find some increasing function $f:\R \to \R$
such that the ball of radius $r$ must contain at least $O(f(r))$
such components.}

\item{Perhaps one should modify this question to measure
generalized Reeb components. For example, the Whitehead
manifold can be obtained as an increasing union of foliated
``plugs'' centered on a single Reeb component.}
\end{enumerate}
\end{rmk}

\begin{qn}
What generalizations of the notion of taut foliation make sense on an
open $3$--manifold?
\end{qn}

\begin{rmk}
This question has different flavors in the topological and the
geometric world.
\begin{enumerate}
\item{For $M$ the interior of a compact manifold with boundary, there
are well--understood generalizations of a taut foliation or a taut
sutured structure; see \cite{Gab83}.}
\item{If $M$ is an infinite cover of some closed $N$, and $N$ admits
a vector field $X$ which is volume--preserving for some volume form on
$N$, a foliation $\F$ of $M$ transverse to the cover $\widehat{X}$ probably
has some useful properties. Maybe one should ask for $\F$ to
be a Lipschitz section of $\widehat{X}$, in the sense that the
transverse flow is uniformly bounded away from $T\F$. In this
case, leaves of $\F$ are approximate minimal surfaces for some
equivariant metric.}
\item{One could weaken the transversality condition on $\widehat{X}$
to an asypmtotic or a statistical condition. For example: we could
ask for every flowline of $\widehat{X}$ to be eventually transverse to
$\F$, Or: assume $\F$ is co--oriented; then one 
could ask for every $\epsilon>0$ there should exist a $t$ such that for
each flowline $\gamma$ of $\widehat{X}$, every
segment of $\gamma$ of length $t$ is transverse to $\F$ in the positive
direction.}
\item{Rather than looking at a vector field $X$, one could ask for a
single loop $\gamma$ in $N$ such that every lift $\widehat{\gamma}$ of
$\gamma$ is transverse to $\F$ in the positive direction, and every
leaf of $\F$ intersects some lift. Again, this could be made into
an asymptotic or statistical condition.}
\end{enumerate}
\end{rmk}

\section{Sublaminations and superlaminations}
\subsection{Essential laminations}

\begin{qn}
Characterize those essential laminations which contain genuine sublaminations.
\end{qn}

\begin{rmk}
Again, {\em characterize} can take several flavors. For example:
\begin{enumerate}
\item{Geometric conditions. For instance, if $\pi_1(M)$ is Gromov--hyperbolic,
a lamination whose universal cover contains a quasi--geodesically embedded
leaf contains a genuine sublamination.}
\item{Topological conditions. For instance, if $\F$ is $\R$--covered or has
one--sided branching, it cannot contain a genuine sublamination.}
\item{Algorithmic conditions. Given a lamination, described somehow in
finite terms, give an algorithm to determine whether it contains a genuine
sublamination.}
\end{enumerate}

Given a branched surface
$\B$, can the algorithm of Agol and Li can be modified to 
tell whether it carries (possibly not fully) a genuine lamination?
(see \cite{Shi97} for related questions). The set of branched surfaces in
a $3$--manifold which fully carry an essential lamination can be constructed, so
the problem reduces to showing, for two branched surfaces $\B,\B'$ whether some
splitting of $\B'$ is carried by $\B$.
\end{rmk}

\subsection{Genuine laminations}

\begin{qn}
Suppose $\Lambda$ is a {\em full genuine lamination}; i.e. it has some complementary region
which is an ideal polygon bundle over a circle. Suppose $M$ is hyperbolic. Is the core
circle of this region isotopic to a geodesic? Does it have a noncoalescable insulator
family?
\end{qn}

\rmk 
\begin{enumerate}
\item{Noncoalescable insulator families are introduced in \cite{Gabcoa} where they
are used to establish (virtual) topological rigidity of hyperbolic manifolds. }
\item{Maybe bundle regions are exactly the wrong kind of region for this analogy. Does
it work better if we ask whether cores of interstitial annuli of non--bundle regions
have insulator families? After all, these regions look more like pared manifolds, where
the annuli are like the pared locus, and therefore maybe have short
geodesic representatives.
Note that such interstitial annuli are unknotted in the universal cover, and therefore,
if $\pi_1(M)$ is residually finite, in some finite cover they lift to
circles isotopic to geodesics. See e.g. \cite{FG01}}
\item{Given some decomposition of $M-\Lambda$ into guts $\mathfrak{G}_i$
and interstitial $I$--bundles, the {\em filling lemma} from \cite{CD01} says that
$\Lambda$ can be {\em filled} to a lamination $\Lambda'$ with gut regions
isotopic to the gut regions of $\Lambda$, and for which every interstitial
$I$--bundle is a product $I \times S^1 \times \R^+$. Thus if $\gamma$ is
the core of a gut region which is a polygon bundle over $S^1$, it is also the
core of a complementary region (of a possibly different lamination) of an ideal
polygon bundle over $S^1$. Perhaps such ideal polygon bundles are {\em artificial};
the more robust examples would be complementary regions to laminations with every
leaf dense --- i.e. minimal laminations.}
\end{enumerate}

\begin{qn}[Thurston]
Suppose $\Lambda$ is a genuine lamination. When can $\Lambda$ be ``filled in'' to
a very full lamination $\Lambda'$? Does it help for $M$ to be hyperbolic?
\end{qn} 

\rmk 
\begin{enumerate}
\item{If $\Lambda'$ is a co--oriented very full lamination,
it can be filled in with monkey saddles to
give a taut foliation. This implies for instance that closed leaves in any
$\Lambda \subset \Lambda'$ should be
Thurston norm minimizing, and $\pi_1(M)$ should admit an action on a universal circle
with no global fixed point. On the other hand, if we do not
require $\Lambda'$ to be co--oriented, there
are no known homological obstructions. There are examples of
homologically trivial closed surfaces which are leaves of very full laminations. 
(\cite{ATun})}
\item{If $\Lambda$ is an essential lamination of $M$, then
the universal cover $(\til{M},\til{\Lambda})$ is topologically a product
$(\R^2,\lambda) \times \R$ where $\lambda$ is a lamination of $\R^2$ (\cite{GKlam}).

Since $\lambda$ is proper, there is a natural circular ordering on the ends of the leaves
of $\lambda$, but this depends strongly on the (non--unique) choice of $\lambda$.
If one can choose an equivariant representation $(\R^2,\lambda)\times \R$ of 
$(\til{M},\til{\Lambda})$, one obtains a
representation of $\pi_1(M)$ in $\homeo(S^1)$ with infinite image. It follows by
\cite{BRW} that $\pi_1(M)$ is isomorphic to a subgroup of $\homeo(S^1)$.
This happens for instance when $\Lambda$ is tight and has solid torus guts.
See \cite{CD01}.}
\end{enumerate}

\subsection{Loosesse laminations}

\begin{defn}
A lamination $\Lambda$ is {\em loosesse} if it satisfies all the properties of
essentiality except possibly the {\em no compressing monogons} condition.
This concept was introduced by Gabai.
\end{defn}

\begin{qn}[Gabai]
Are loosesse laminations good for anything?
Are leaves of the universal cover of a loosesse lamination properly
embedded? If $M$ contains a loosesse lamination, does $\til{M} = \R^3$?
\end{qn}

\rmk Brittenham has observed that for {\em branched surfaces}, the ``no compressing
monogons'' condition can be substantially weakened to ``no monogon bundles over $S^1$''.

If $M$ contains a loosesse lamination $\Lambda$ for which $\til{M} = \R^3$, is there an
analogue of \cite{GKlam} to the effect that $(\til{M},\til{\Lambda})$ is topologically
$(\R^2,\lambda) \times \R$ for some (not necessarily proper) lamination $\lambda$ of
$\R^2$?

\begin{qn}
Give an example of a lamination in an atoroidal manifold
--- perhaps loosesse --- which can never be realized
by minimal surfaces for any metric, but which certifies some useful topological
property of $M$.
\end{qn}

\rmk It is unknown whether every essential lamination can be realized by minimal
surfaces for some metric; in some sense they seem to function as though they can be.
But laminations containing compressing monogons can {\em never} be minimal surfaces
for any metric.

\section{Branched surfaces and triangulations}
\begin{qn}[Agol]
Characterize branched surfaces embedded in $3$--manifolds which can be
non--trivially split to a homeomorphic copy of themselves.
\end{qn}

\begin{rmk}
\begin{enumerate}
\item{If $\B$ is such a branched surface, one can set $\B_0 = \B$ and
$\B_i$ the result of performing such a splitting on $\B_{i-1}$. If one
fixes a triangulation $\tau$ transverse to $\B$, this procedure gives
rise to a partial lamination $\L$ carried by $\B$ (``partial'' because
the injectivity radius of the subsurface which is being split open might
be bounded above). If $\B$ carries a lamination, is
$\L$ contained in some lamination $\Lambda$
carried by $\B$?}
\item{If $\B$ supports a nontrivial transverse measure $\mu$, it
supports a nontrivial measure projectively invariant under splitting
(this follows from the fact that the set of nontrivial projective
measures is a convex set, and the splitting operation induces a
projective transformation of such a set, which necessarily has a
fixed point). One can also look at nontrivial transverse measures
twisted by some fixed holonomy representation of the fundamental group of $\B$.
Is there some underlying geometric meaning of such measures? A basic example
is a branched surface carrying the stable or unstable laminations of
a pseudo--Anosov flow transverse to a surface bundle (such a branched
surface can be derived from a sequence of splittings and collapses of an
supporting train track).}
\end{enumerate}
\end{rmk}

\begin{qn}[Agol]
Develop a theory of hierarchies for branched surfaces.
\end{qn}
\begin{rmk}
A branched surface $\B$ would be decomposed sequentially along train tracks
into simple pieces. The idea would be to try to prove a gluing theorem about
how essential laminations carried by the simple pieces of the hierarchy
glue together to give essential laminations in the branched surface.

Something like this is done by Li in \cite{Li00} for branched surfaces without
sink disks.
\end{rmk}

Recall the notion of a {\em taut ideal triangulation} as introduced by Lackenby
in \cite{Lacktit},
as a combinatorial refinement of certain ideal triangulations of hyperbolic
$3$--manifolds with torus cusps. Such objects give rise to branched surfaces
which are very simple --- their worst singularities look like train tracks $\times I$.

\begin{qn}[Dunfield]
Which boundary slopes are realized by essential laminations carried by a fixed
taut ideal triangulation? Give an algorithm.
\end{qn}

\rmk Though these branched surfaces do not have sink disks, they have many
{\em half sink--disks}, so not all of Li's constructions for manifolds with
boundary can be made to work. These are good candidates for the ``automatic
laminations'' in question~\ref{automatic_lamination}.

\begin{qn}[Schleimer]
When does a Haken sum operation make sense for a pair of laminations
in normal form with respect to a fixed triangulation?
\end{qn}

\rmk For transversely measured laminations, this question is less interesting.
Although, see Hatcher \cite{Hatmea}.

A prerequisite must be that the laminations intersect each tetrahedron
in the same kind of quadrilateral pieces. This is not too far from the
condition that they are simultaneously carried by a fixed branched surface.
So a better question might be --- what is the
structure of the set of laminations carried by a given branched surface (up to
monotone equivalence, say)? Is there a natural topology for which it is star--like?
Or contractible?

\subsection{Normal surfaces}

The theory of normal surfaces was introduced by Haken in \cite{Hak61}, although
similar ideas were introduced in a more limited scope by Kneser. Any foliation
can be normalized relative to some triangulation \cite{Ben97}.

\begin{qn}
Let $M$ be a $3$--manifold, and $\Lambda$ an essential lamination.
Let $C$ be a cycle representing the fundamental class of $M$. Is there a cycle
$C'$ with the same Gromov norm as $C$, and another essential lamination $\Lambda'$ which
is normal with respect to $C'$?
\end{qn}
\begin{rmk}
\begin{enumerate}
\item{If $\Lambda$ is an incompressible surface, $C$ can be normalized by \cite{Agol}. In fact,
Agol points out in his paper that the same method works if $\Lambda$ is merely tight.
Presumably, the cycle $C'$ and the lamination $\Lambda'$ will be derived from $C$ and
$\Lambda$ by some process of straightening the simplices of $C$ relative to $\Lambda$,
and compressing and evacuating (\cite{Gab99})
$\Lambda$ with respect to its intersection with $C$.

A positive answer would extend Agol's volume bounds to manifolds with essential (not
necessarily tight) laminations.}
\item{Does it help to assume all the coefficients of $C$ are positive? In this case,
for any singular simplex $\sigma:\Delta^3 \to M$ in the support of $C$, and any leaf
$\lambda$ of $\Lambda$ intersecting $\sigma(\Delta^3)$, we can {\em develop} (not uniquely)
$C$ along any path $\gamma \subset \lambda$ which does not pass through any edges or vertices.
For, once $\gamma$ comes to a face $\tau$ of $\sigma(\Delta^3)$, there is another singular
simplex $\sigma':\Delta^3 \to M$ whose boundary absorbs some of the mass of $\tau$; then
$\Lambda$ can be developed along $\gamma$ into $\sigma'$, and so on inductively. This
gives a particular order in which to try to compress and evacuate the leaves of $\Lambda$
relative to $C$.}
\item{If $C$ is a virtually embedded cycle --- i.e. it comes from a triangulation in some
finite cover $\widehat{M}$, $\Lambda$ can be lifted to $\widehat{\Lambda} \subset \widehat{M}$
and evacuated there. Can it be evacuated equivariantly?}
\end{enumerate}
\end{rmk}

\begin{defn}
A {\em taut local orientation} is a choice of ordering for the vertices
of each tetrahedron in a triangulation with the following properties:
\begin{enumerate}
\item{The star of each vertex
is ordered compatibly with a local foliation in normal form}
\item{Every oriented loop is homotopically essential}
\end{enumerate}
\end{defn}

Such objects are dual to branched surfaces with certain properties; in particular,
certain examples give rise to branched surfaces which carry nothing,
but fully carry essential laminations in finite covers. See \cite{Cal00a} or \cite{Calp2}.

\begin{qn}
Suppose $\B$ is a branched surface in $M$ which is dual to a taut local orientation.
Is there a finite cover of $M$ in which the pullback of $\B$ fully carries a
lamination? What about an amenable cover?
\end{qn}

Say an {\em automatic lamination} is given by the following data:
\begin{itemize}
\item{A branched surface $\B$ with some fixed cell decomposition, so that
there is a notion of combinatorial paths supported by this surface
which can be enumerated.}
\item{A lamination $\Lambda$ fully carried by $\B$.}
\item{A finite state automaton $\sf A$ which recognizes combinatorial
paths in $\B$ which represent geodesic paths in leaves of $\Lambda$.}
\end{itemize}

\rmk Note that if $\Lambda$ is essential and $M$ is atoroidal, leaves of
$\Lambda$ are uniformly $\delta$--negatively curved, so such an automaton
exists relative to an oracle which recognizes combinatorial paths in
$\B$ which are carried by some leaf of $\Lambda$.

\begin{qn}\label{automatic_lamination}
Do branched surfaces without sink disks carry automatic laminations?
\end{qn}

\begin{rmk}
\begin{enumerate}
\item{This seems possible, insofar as Li's construction of a lamination carried by
such a branched surface is local and explicit, and describes a procedure for
recursively splitting open the branched surface. One subtle issue might be
the choice of an extension over foliated circle bundles over punctured surfaces of
positive genus: the existence of such an extension follows from the fact that
every homeomorphism of $S^1$ is a commutator of length at most $2$.
But it is not clear whether a
recursive description of such a commutator can be given ``automatically''.}
\item{Every time a combinatorial path passes over a source branch, there are two
possible routes it can take. The sequence of splittings (positive, negative)
so obtained is like a continued fraction expansion for an irrational number.
The automaton $\sf A$ can only remember a bounded length of this
sequence; this constrains the possible transverse structures supported
by an automatic lamination. In particular, any automatic lamination should admit a natural
transverse bilipshitz structure.}
\item{This is related to question~\ref{automatic_order} on the
existence of {\em automatic left--orderings} on
finitely presented groups.}
\end{enumerate}
\end{rmk}

\begin{qn}
Give a useful definition of {\em thin position} for an
embedded graph $\Gamma \subset M$ with respect to a taut foliation $\F$.
If $\Gamma$ is the $1$--skeleton $\tau^1$ for a triangulation, can one find
an isotopy such that the leaves of $\F$ are made up of polyhedral disks of bounded
index?
\end{qn}

\begin{rmk}
For $\F$ a foliation by closed surfaces, one can make the leaves
either normal or almost normal. If $\F$ is finite depth,
one can normalize the depth $0$ leaves. The depth $1$
leaves are asymptotic to the depth $0$ leaves, so they may be assumed to
be normal away from a compact piece; some kind of relative minimal sweepout
argument might produce an isotopy so that these leaves are almost normal.
Continuing inductively, it seems that one can bound the index of the disks
appearing in a leaf by a function of their depth. This argument should be
made more precise.
\end{rmk}

\section{Leaf spaces and transverse structures}

The most significant nonexistence results in this area are the main theorems of
\cite{RSS} and \cite{Fen_nolam} which give infinitely many examples of hyperbolic
$3$--manifolds $M$
for which there is no nontrivial
action of $\pi_1(M)$ on a simply--connected $1$--manifold, respectively
order tree.

\begin{qn}[Thurston]
Suppose $M$ is irreducible. Suppose further that
$\pi_1(M)$ admits a nontrivial action on $\R$. When does $M$ admit a taut foliation
with a transverse $(\pi_1(M),\R)$ structure?
\end{qn}
\begin{rmk}
\begin{enumerate}
\item{Given a representation $\rho:\pi_1(M) \to \homeo(\R)$, let $(E_\rho,\F_\rho) \to M$
be the associated foliated bundle. A section $s:M \to E_\rho$ transverse to $\F_\rho$
pulls back $\F_\rho$ to a foliation of $M$ with a transverse $(\pi_1(M),\R)$ structure.}
\item{Given an action of a group $\Gamma$ on a simply connected $1$--manifold $L$,
when is there an immersion $i:L \to \R$ such that the action of $\Gamma$ descends to
an action on $\R$? One should allow some flexibility in addressing this question, so that
for instance one should be allowed to substitute monotone equivalent actions. Does it
make a difference to weaken the quality of $i$ so that it is not necessarily an immersion,
but merely (partial) order preserving? Under what circumstances is there a {\em unique}
maximally ordered $\Gamma$--invariant order--preserving
quotient $L \to L_o$? What if one imposes extra
conditions on $\Gamma$, for instance that it should be amenable (especially if one allows
monotone equivalent actions)?}
\item{One can obviously
generalize this question to actions on order trees.}
\end{enumerate}
\end{rmk}

\begin{qn}
Suppose $\F$ is an $\R$--covered foliation of an atoroidal $3$--manifold $M$. Is the
holonomy representation $\rho_H$ of $\pi_1(M)$ on $\R$ conjugate to a group of
coarse $1$--quasi--isometries? i.e. is there a positive valued
function $C:\pi_1(M) \to \R$ such that
$$| d(\rho_H(\alpha)(p), \rho_H(\alpha)(q)) - d(p,q) | \le C(\alpha)$$
for all $\alpha \in \pi_1(M)$ and $p,q \in \R$?
\end{qn}

\rmk 
\begin{enumerate}
\item{If the value of $C$ can be chosen independent of $\alpha$, it follows that
any two leaves of $\til{\F}$ are contained within bounded neighborhoods of each other.
Such foliations are said to be {\em uniform}. All known examples of $\R$--covered
but nonuniform foliations of atoroidal $3$--manifolds satisfy the condition above.}
\item{If $M$ is allowed to be toroidal, the action of $\pi_1(M)$ on $L$ can fail to
satisfy the condition above. It is compatible with all known examples that for
toroidal manifolds, the action of $\pi_1(M)$ should always be
conjugate to a group of coarse quasi--isometries.}
\end{enumerate}
See \cite{Thu97} and \cite{Cal99} for more details.

\begin{qn}
Suppose $M$ is atoroidal and admits a taut foliation. Must it admit an $\R$--covered
foliation?
\end{qn}

\rmk
\begin{enumerate}
\item{There are many examples in \cite{CD01} of small volume hyperbolic manifolds which
admit essential laminations, but which do not admit $\R$--covered foliations. In
fact, their fundamental groups admit no nontrivial representation in $\homeo(\R)$.}
\item{If one is prepared to pass to a finite cover, one probably cannot
find an example based on the nonexistence of group
actions on certain $1$--manifolds. For, with suitable homological
hypotheses (which can be achieved by passing to a finite cover), the fundamental
group of a tautly foliated manifold acts on $\R$ without a global fixed
point. See \cite{CD01} for details.
On the other hand, \cite{RSS} gives many examples of fundamental groups
of hyperbolic $3$--manifolds which admit no actions on $\R$ without a global fixed
point, but which admit taut foliations and therefore have faithful circle actions.
The Euler class of these circle actions are torsion. In general, for any $3$--manifold
$M$ and for any $2$--plane field $\xi$ on $M$ there is a $3$--manifold $N$ together with a
map $f:N \to M$ that induces an isomorphism on homology (even over $\pi_1(M)$) such that
$f$ is transverse to $\tau$, and such that the induced $2$--plane field $f^*\xi$ on
$N$ is homotopic to a taut foliation. For, Thurston's construction \cite{Thu73} lets
us construct a foliation $\F$ on $M$ in a neighborhood of the $2$--skeleton of an
anti--oriented triangulation transverse to $\tau$. By drilling out spiralling paths, one can
arrange for the complementary regions to the foliated region to be solid tori $S_i$. Then
$N$ can be obtained from $M$ by drilling out these solid tori $S_i$ and filling in circle
bundles $B_i$ over punctured surfaces. Any homeomorphism of $S^1$ can be realized as
a commutator of bounded length, so the foliations of $\F|_{\partial S_i}$ can be
extended to foliations of $B_i$ transverse to the circle fibers. This gives a
demonstrably taut foliation in the appropriate homology class on a manifold with the
required properties. This construction is due to Thurston \cite{Thurmk}, and
gives many examples of circle actions of fundamental groups of
hyperbolic manifolds whose Euler class is torsion but non--zero.}
\item{Brittenham, Naimi and Roberts (\cite{BNR}) have produced examples of graph
manifolds which admit taut foliations but do not admit $\R$--covered foliations.}
\item{A related question from \cite{Gabprob} is: if $M$ admits an
essential lamination, does it admit a tight essential lamination? (an
essential lamination is {\em tight} if the leaf space of the universal
cover is Hausdorff). Quasigeodesic laminations in hyperbolic
manifolds are tight, so a positive answer to question~\ref{quasigeodesic}
would give the most positive answer to this question.}
\end{enumerate}

\begin{qn}[Agol]\label{SL2R_laminations}
For a fixed manifold $M$, describe the structure of the set of all essential
laminations with a transverse $\til{SL(2,\R)}$ structure.
\end{qn}
\begin{rmk}
\begin{enumerate}
\item{For a fixed triangulation $\tau$, the set of normal foliations supported by
$\tau$ with a $\til{SL(2,\R)}$ structure can be algorithmically described, and
has the natural structure of a semi--algebraic set. By Gabai \cite{Gab99}, there
is a fixed triangulation $\tau$ of $M$ so that every essential lamination of $M$
can be made normal with respect to $\tau$; in fact, there is a constructive
procedure to find a finite set of branched surfaces which fully carry every nowhere
dense essential lamination in $M$.}
\item{One can study more generally the space of actions of $\pi_1(M)$ on
$SL(2,\R)$--trees, which one could ask to be either partially ordered or not.
In this context, it is probably useful to differentiate
between {\em complete} and {\em incomplete} trees. Any properly
embedded interval $I$ in an $SL(2,\R)$--tree, has a well--defined {\em coarse length},
a non--negative integer $l(I)$ which is the number of times $I$ wraps around
$S^1$ under projection. An $SL(2,\R)$--tree is
{\em complete} if every properly embedded tight ray $r$ has infinite
coarse length. In particular, if some element $\alpha$ stabilizes some
point $p \in r$, it stabilizes
infinitely many points on $r$ which exit the noncompact end of $r$. This problem is
probably a good warm--up to understanding the space of actions of $\pi_1(M)$ on
arbitrary order trees or $1$--manifolds.}
\item{Suppose $\Gamma$ acts on an $\R$ order tree $T$ in such a way that
for every nontrivial $\gamma$, and for every properly embedded copy $l$ of
$\R$ in $T$ invariant under $\gamma$, the set of fixed points of $\gamma$
is isolated, and is either empty or exits both noncompact ends of $l$,
and the translation direction of $\gamma$ on the complementary intervals
of $l$ alternate. Suppose furthermore that there is {\em no} sequence $\gamma_i$ with
consecutive fixed points $p_i,q_i$ contained in a fixed segment $I$ for which
$p_i,q_i \to r$ for some $r \in I$. Does $T$ carry the natural structure
of an $SL(2,\R)$--tree?}
\end{enumerate}
\end{rmk}

\begin{qn}
Suppose $M$ admits a minimal taut foliation. What is the best analytic (transverse) quality
of a taut foliation it admits? Can we find a minimal foliation such that the holonomy
groupoid is of type $\text{\rm III}_\lambda$ for some algebraic $\lambda$?
What about if one asks for a foliation monotone equivalent to the first? Homotopic?
\end{qn}
\rmk 
There are well--known topological obstructions to smoothing transverse structures.
\cite{Pix77} gives some important examples of foliations which can be made transversely
$C^1$ but not $C^2$. At the level of Haefliger structures, the obstructions are
related to homology of various groups of homeomorphisms of $\R^n$
(\cite{Hae57},\cite{Hae62},\cite{Tsu90}).

A pseudogroup $\G$ of transformations of a measure space $(X,\mu)$ is said to be
of {\em type $\text{\rm III}_\lambda$} if for each $g \in \G$, the
Radon--Nikodym derivative $d(g_*\mu)/d\mu$ is well--defined on a set of full
measure and takes values in $\langle \lambda \rangle$, the abelian group of
powers of $\lambda$. Such measures and actions arise frequently from
automatic structures; the theory of Patterson--Sullivan measures is also relevant.

An obvious obstruction is that $\pi_1(M)$ must act on {\em some} $1$--manifold
with this quality. Question~\ref{automatic_lamination} and
question~\ref{universal_quality}
are relevant here. See \cite{Con94} for an abstract discussion of transverse measure
theory for foliated manifolds. 

\begin{qn}
Is there a universal constant $c$ such that a hyperbolic $3$--manifold $M$
whose fundamental group $\pi_1(M)$ can be ordered out to radius $c$ can be
left--ordered? Or weaker, is there an effective method to compute such a
$c(M)$ for a given $M$?
\end{qn}
\begin{rmk}
\begin{enumerate}
\item{A $3$--manifold group $\pi_1(M)$ is left orderable iff it admits a
faithful representation in $\homeo^+(\R)$. If $M$ admits an $\R$--covered foliation,
the representation $\pi_1(M) \to \homeo(L)$ may not be faithful, but the
kernel is a surface group, and is itself left--orderable, so the representation
is monotone equivalent to a faithful one.}
\item{The set of finitely presented groups which are not left--orderable is
recursively enumerable; on the other hand, the set of finitely presented groups
which {\em are} left--orderable is not. Does hyperbolic geometry, or more
generally word--hyperbolic geometry make a difference here? A word--hyperbolic
group looks more or less free on a large scale; free groups are certainly
left--orderable. A positive answer to this question would give an algorithm
to detect left--orderability (of hyperbolic $3$--manifold groups).
In fact, any $M$ which admits a taut foliation,
a pseudo--Anosov flow, or certain kinds of genuine laminations,
has a faithful representation $\pi_1(M) \to \homeo(S^1)$ coming from a universal
circle; if $M$ is a rational homology sphere, this representation lifts to
$\til{\homeo(S^1)}$ when restricted to the commutator subgroup, a finite index
subgroup. Such periodic representations might be easier to detect or rule out
than arbitrary ones. In particular, in \cite{CD01} it is proved that the
Weeks manifold does not admit a pseudo--Anosov flow, or a tight essential lamination.
On the other hand, there are examples of $3$--manifolds which act faithfully on $S^1$
but not on $\R$. See \cite{CD01} for more details.}
\end{enumerate}
\end{rmk}

\begin{qn}\label{automatic_order}
Let $\sf{T}$ be some class of abstract computers; e.g. finite state automata,
Turing machines, Turing machines relative to some oracle $O$, etc.
A {\em $\sf{T}$--order} on a group $G$ is a left--invariant order such that
there is a machine $T \in \sf{T}$ which recognizes the positive cone $G^+ \subset G$.
What kinds of $\sf{T}$--orders are possible for fundamental groups $G$ of
hyperbolic $3$--manifolds? 
\end{qn}
\begin{rmk}
\begin{enumerate}
\item{A special case consists of the class $\sf{T}$ of finite state automata. One
can ask whether there exists an automatic order on $G$; i.e. whether there is an
FSA which can detect {\em geodesic} words with values in the positive cone. Of course,
$G$ must already be an automatic group for there to be any chance at all. It turns out
that this class of groups is very small. There are automatic orders on $\Z^n$ for all
positive integers $n$, but there is no automatic order on $F_2$, the free group on
two generators. A proof is as follows. Let $F_2 = \langle a,b \rangle$ be
a generating set. A left--invariant order on $F_2$ is more or less equivalent to
a faithful action of $F_2$ on $\R$ by orientation--preserving homeomorphisms.
We may also assume this action is minimal. Let $0$ be a point with trivial stabilizer,
so that the positive cone are the elements $g$ with $g(0)>0$.

Suppose $w$ acts as a positive translation on a $w$--invariant interval $I$, with
lowest point $p<0$ and highest point $p'>0$. Since the action is minimal, there is
$q \in I$ and a word $v$ such that $v(q) = p'$. Let $x_1$ be a word with $x_1(0)=p_1^+$,
where $p_1^+>p$ is very close to $p$. Then 
$$w^{-n}vw^nx_1(0)>0$$ for $n>N_1$, and
$$w^{-n}vw^nx_1(0)<0$$ 
for $n \le N_1$ for some very large $N_1$.
We can choose $x_i$ with $x_i(0) \to p$ for which there
are corresponding $N_i \to \infty$. If $\sf{A}$ is a fixed finite state automaton,
eventually $N_i$ exceeds the memory capacity of $\sf{A}$, and $\sf{A}$ is
incapable of distinguishing $w^{-N_i}vw^{N_i}x_i$ from $w^{-N_i-1}vw^{N_i+1}x_i$.

Conversely, if
every $w$ fixes at most $1$ point on $\R$, the image of the group is solvable,
and therefore not faithful.

Since every hyperbolic $3$--manifold group contains a quasigeodesically embedded
$F_2$, this implies that no hyperbolic $3$--manifold group admits an
automatic order.}
\item{If a finitely generated automatic 
group $G$ is left--orderable, given a lexicographic ordering on representatives
of $G$ there is a well--defined {\em lexicographically first} left--order.
What is the algorithmic quality of this order? Note that given an automatic
group $G$, there is an algorithm which either certifies that $G$ is not
left--orderable, or recursively enumerates the words in the lexicographically
first left--order.}
\item{The general question of whether or not a finitely presented group is
left--orderable is undecidable. But even if a group is left--orderable, it
might be much harder to describe an explicit left--order than merely to
certify orderability.}
\end{enumerate}
\end{rmk}

\subsection{Universal circles}

\begin{qn}\label{lam_question}
Let $\Lambda^\pm$ be a pair of laminations of $S^1$ which are transverse to each
other and have finite area complementary domains.
Suppose $\Gamma$ is a group of automorphisms of $S^1$
which preserves $\Lambda^\pm$ and acts minimally on
the leaves of either lamination. When is
$\Gamma$ commensurable with $\pi_1(M)$ for $M$ a hyperbolic $3$--manifold?
\end{qn}
\begin{rmk}
\begin{enumerate}
\item{If $\Gamma = \pi_1(M)$ for some $M$, must $M$
admit a pseudo--Anosov flow $X$ for which
the associated action on the universal circle of $X$ is conjugate to the given action of
$\pi_1(M)$ on $S^1$? Actions of $3$--manifold groups
on circles associated to pseudo--Anosov
flows are constructed in \cite{CD01}.}
\item{What properties should a lamination $\Lambda$ or collection of laminations
of $S^1$ satisfy to ensure good geometric control over the subgroup of $\homeo(S^1)$
preserving $\Lambda$? For instance, let $\Lambda^\pm$ be as above and $\Gamma$ a
group of automorphisms of $\Lambda^\pm$: is $\Gamma$ word hyperbolic?  
Is there a natural topology on the space of such pairs
$\Lambda^\pm$ for which a large group $\Gamma$ of automorphisms exists?
What do non--trivial convergent sequences
(of commensurability classes of hyperbolic
$3$--manifolds) look like in such a topology?}
\item{Given a pair of laminations $\Lambda^\pm$ of $S^1$ as above, the natural
quotient of $S^1$ by the leaf relations of $\Lambda^\pm$ is topologically $S^2$.
If $\Gamma$ is a group of automorphisms of $\Lambda^\pm$, we get a representation
$\Gamma \to \homeo(S^2)$ and a flat bundle over a $B\Gamma$.
The topological group $\homeo(S^2)$ is homotopy equivalent
to $O(3,\R)$ by Smale; are the cohomology classes
pulled back via the associated classifying map $B\Gamma \to BO(3,\R)$ bounded in
norm in terms of combinatorial data that can be read off from $\Lambda^\pm$?}
\end{enumerate}
\end{rmk}

\begin{qn}
What possibilities are there for universal circles $S^1_\u$ for a fixed manifold?
For a fixed foliation? For what taut
foliations is there a {\em unique minimal} universal circle?
\end{qn}

\begin{rmk}
\begin{enumerate}
\item{A universal circle $S^1_\u$
for a taut foliation $\F$ formally gives rise to a
nontrivial invariant lamination of $S^1_\u$
for each direction in which $\til{\F}$ branches (\cite{Calp3}).
Given an action of a group $G$
on a circle $S^1$ preserving a lamination $\Lambda$, one can modify the action
in quite drastic ways: one can take an orbit class of complementary region, and
``flip'' the lamination along every edge of every region in the orbit. If $T$ is
the dual planar tree to $\Lambda$, this is equivalent to re--embedding $T$
in the plane in such a way that the circular order is reversed at every vertex
in an orbit class (see \cite{CD01}).
If there is more than one orbit class of complementary region, this
operation is not a topological conjugacy, and the circle actions in question
are quite distinct. One can almost certainly use this operation to produce examples
of minimal circle actions of hyperbolic $3$--manifold groups
which are not universal circles for any foliations.}
\item{There is a unique minimal universal circle for $\F$ an $\R$--covered foliation,
or one with one--sided branching. These minimal universal circles are rigid under
certain kinds of deformations of the foliation, in many cases.}
\item{The best answer would do two things: firstly it would give a good description of
the ``representation space'' of $\pi_1(M)$ in $\homeo(S^1)$ and secondly
it would explicitly characterize the subspace of representations of geometric origin.

It is probably unreasonable to expect an algebraic structure on the representation
space, but a {\em polyhedral structure} might be more promising, by analogy with
question~\ref{polyhedral_pseudos}.  A sequence of isotopies of a fixed
pseudo--Anosov flow degenerating to a flow corresponding to a higher codimension face
should give rise to a monotone relation between the universal circles of the leaf
space of the flows.}
\item{Universal circles for certain kinds of genuine laminations
and for pseudo--Anosov flows are constructed in \cite{CD01}. Perhaps every
essential lamination gives rise to a universal circle.}
\end{enumerate}
\end{rmk}

\begin{qn}\label{universal_quality}
What is the best analytic quality for the action of $\pi_1(M)$ on
a universal circle $S^1_\u$?
\end{qn}
\begin{rmk}
\begin{enumerate}
\item{It is probably unreasonable to expect the action to be $C^2$ if $M$ is
hyperbolic. On the other hand, {\em any} minimal action of a finitely
generated group on $S^1$ is conjugate to a Lipschitz action, with respect
to a harmonic measure. For an introduction to the rich panorama of
subgroups of homeomorphisms of $S^1$ stretching between Lipschitz and
$C^2$, see \cite{Tsu95}.}
\item{The following related question was suggested by Barry Mazur.
Let $G$ be a finitely presented group, and let $X,Y$ be compact smooth manifolds with
$\text{dim}(X) < \text{dim}(Y)$. Suppose
$\rho_X:G \to \homeo(X)$ and $\rho_Y:G \to \homeo(Y)$ are minimal actions such
that there is a surjection $f:X \to Y$ which intertwines $\rho_X$ and $\rho_Y$.
When are $\rho_X,\rho_Y$ actions by diffeomorphisms for some appropriate smooth structures
on $X,Y$, and how can one find construct such a smooth structure on $X$ from a smooth
structure on $Y$ and from $f$? This question is most interesting when $X$ is a circle
or an interval, but there are probably many interesting examples of $Y$. What seems
plausible is that high dimensional
cohomology of $G$ pulled back from the action on $Y$ represents an obstruction to
stiffening the analytic quality of the action on $X$.}
\end{enumerate}
\end{rmk}

\section{Classical $3$--manifold theory}

\begin{qn}
Is there a {\em universal} transverse surgery description of tautly
foliated manifolds, in the
sense that there is a fixed $M$ such that for every tautly foliated manifold $N,\F$
there is a link $L \subset N$ transverse to $\F$ so that $M$ is obtained from $N$ by
surgery on $L$?
\end{qn}
\begin{rmk}
\begin{enumerate}
\item{The idea is to get an analogue of Lickorish's theorem
(that every closed $3$--manifold
is surgery on a link in $S^3$) in the foliated context. One should also have in
mind Thurston's observation (see \cite{Ada99}) that the process of repeatedly
drilling out the shortest geodesics in a sequence of hyperbolic
manifolds tends to stabilize to a ``universal'' sequence.
Question~\ref{short_geodesics} is probably relevant
to an approach along these lines. One should concentrate on taut foliations for this
question, so that there are enough transverse knots to ``see'' the whole topology
of $N$. One idea might be to start with a triangulation $\tau$ in normal
form with respect to $\F$, and look at the interaction of the
associated Heegaard decomposition with $\F$.}
\item{We can modify the question to ask whether, given two tautly foliated
manifolds $M_1,\F_1$ and $M_2,\F_2$, there are links $L_i$ transverse to
$\F_i$ such that $M_1\backslash L_1 = M_2\backslash L_2$. }
\end{enumerate}
\end{rmk}

\begin{qn}
Give a collection of fundamental operations on foliations and an explicit family of
{\em base foliations} such that every tautly foliated manifold $M,\F$ is obtained
from one of the base family by repeated application of fundamental operations.
\end{qn}

\begin{rmk}
\begin{enumerate}
\item{It is unclear if even the following naive list of operations

\begin{itemize}
\item{branched covers and pushforward}
\item{cut--and--shear}
\item{torus connect sum}
\item{tangential surgery}
\item{surface bundle plumbing along transverse train--tracks}
\item{monotone deformations}
\item{cut--and--paste furrows along loops with linear holonomy}
\end{itemize}

and base foliations

\begin{itemize}
\item{foliations of circle bundles}
\item{transversely measured foliations}
\item{finite depth foliations and those representing irrational rays in $H_2$}
\item{foliations derived from alternating knots by Delman--Roberts methods}
\end{itemize}

would fail to suffice. The idea is really to develop some more robust
invariants of foliations which vary either not at all, or in some well--understood way
under basic operations, and at the same time to give a procedure for recognizing
{\em orphans} with respect to a given list of operations.}
\item{A model theorem along these lines is the recent result of Noah Goodman and,
independently, Giroux, improving work of Gabai, which says that the set of
fibered links in a homology $3$--sphere $M$ are related by a sequence of
elementary stabilizations and destabilizations, where the stabilization
operation is plumbing with a Hopf band.}
\end{enumerate}
\end{rmk}

\subsection{Persistence questions}

\begin{qn}
What is the most general class of knots to which the techniques of Delman--Roberts
(in constructing persistent laminations) can be extended?
\end{qn}

\begin{rmk}
The main reference is \cite{DRprop}.

The idea of this question is that the lamination should
be constructed in an essentially algorithmic
way from a {\em projection} of the knot. A related question asks
what information essential laminations can give about projections
of knots; for instance, when can foliated techniques be used to
certify that a knot projection has minimal number of crossings?
Is the problem of minimal crossing number for knot projection NP--complete?
\end{rmk}

\begin{qn}[Delman]\label{delman}
Suppose $K$ is a non--torus alternating knot. Then essential laminations can be
constructed which realize every (nontrivial) boundary slope. Can essential
laminations be constructed with an even sided bundle complementary region
containing $K$, so that every nontrivial surgery can be filled in with a monkey
saddle?
\end{qn}

\subsection{Distance of Heegaard splittings}

\begin{qn}
It is known \cite{Hart} that if a $3$--manifold $M$ contains an
essential surface of genus $g$, the {\em distance} of any Heegaard
splitting of $M$ has distance at most $2g$. Does the ``distance
filtration'' put any useful structure on the essential laminations supported
by a given $M$? i.e. if $M$ admits Heegaard splittings of distance at
least $2g$, what can one say about the essential laminations $\Lambda$
contained in $M$? 
\end{qn}
\begin{rmk}
A preliminary technical issue would be to decide what a {\em tight} intersection
of $\Lambda$ with the splitting surface $\Sigma$ should be, and then try to
tighten any given $\Lambda$ with respect to $\Sigma$, perhaps by means of an
infinite process, by analogy with Brittenham's principle (see \cite{Bri95} and
\cite{Gab99}). 
\end{rmk}

\section{Hyperbolic geometry}
\subsection{Asymptotic geometry of leaves}

\begin{qn}\label{asymptotically_separated}
Suppose $\F$ is a taut foliation of a hyperbolic $3$--manifold $M$ with
two--sided branching.
Must there be a leaf $\lambda$ of $\til{\F}$ whose complement contains
an open halfspace of $\H^3$ on either side? We call such a leaf {\em asymptotically
separated}.
\end{qn}

\rmk This question is equivalent to the condition that the limit set of
some leaf of $\til{\F}$ is not all of $S^2_\infty$. A conjecture of Fenley
says that $\lambda_\infty = S^2_\infty$ for some leaf $\lambda$ of $\til{\F}$
iff $\F$ is $\R$--covered.

\begin{qn}
Do leaves of $\til{\Lambda}$ for $\Lambda$ an essential lamination have the
{\em continuous extension property?} More generally, what is the relationship
between the action of $\pi_1(M)$ on various ideal boundaries of $\til{M}$
arising from the foliated structure (e.g. universal circles) and the ideal
boundaries arising from the geometry of $\til{M}$.
\end{qn}

\begin{rmk}
\begin{enumerate}
\item{The seminal result in this area is the theorem of Cannon and Thurston \cite{CT}
that leaves of surface bundles over $S^1$ extend continuously to sphere filling
curves. The method of Cannon and Thurston has been greatly extended and improved
by Fenley to deal with many (most?) finite depth foliations, and certain other
taut foliations \cite{Fen92},\cite{Fen98},\cite{Fen99}.}
\item{A {\em measurable} extension property is much easier to establish, and
follows from exponential growth of hyperbolic spaces, and the fact that a leaf
of an essential lamination in the universal cover is uniformly properly embedded in
its $\epsilon$ neighborhood for some $\epsilon$. See \cite{Calmea}.
Note that even in cases
where the continuous extension property is known leafwise (e.g. finite depth
foliations) it is not known for all or any universal circles (although Fenley
has some unpublished partial results on this question).}
\end{enumerate}
\end{rmk}

\begin{qn}\label{dimension_depth}
Suppose $\F$ is a finite depth foliation of a hyperbolic $3$--manifold.
What is the relationship (if any) between the Hausdorff dimension of the limit set
of a leaf $\lambda$ of $\til{\F}$ and the depth of $\F$ or $\lambda$?
\end{qn}

\begin{rmk}
\begin{enumerate}
\item{Conjecturally, a foliation is $\R$--covered iff the limit set of any
(and therefore every) leaf is all of $S^2_\infty$. A proper foliation of depth $>0$
contains a quasigeodesic leaf; in particular, depth 0 is distinguished from
any other finite depth by the property that the limit set of every leaf has
Hausdorff dimension 2. More generally, if $M' \subset M$ is the union of the leaves
of highest depth, the components of $\til{M}'_\infty$ and $\lambda_\infty$ for
leaves $\lambda$ of highest depth should agree.}
\item{For depth one leaves, the components
of $M'$ are the interiors of compact hyperbolic manifolds with quasigeodesic boundary.
For such manifolds, elements of $\pi_1$ can be enumerated by a finite state automaton.
Does this mean that the Hausdorff dimension of the limit set is
a special value of an $L$--function?}
\item{One can numerically estimate the Hausdorff dimension of the limit set
of a noncompact leaf, using a harmonic measure for the universal circle
to determine which directions
at infinity (in the intrinsic geometry of the leaf) are most distorted in the
ambient geometry. This should give
a fast numerical algorithm to tell whether or not a foliation of a hyperbolic
manifold is of depth one.
On the other hand, when $M$ has cusps, there are numerical problems with this approach,
since it is slow to enumerate elements of $\pi_1(M)$ which exit a cusp.
This is especially troubling, since the case one would be most interested in is
foliations of knot or link complements. It is an interesting problem to find therefore
a fast numerical algorithm to estimate the Hausdorff dimension of the limit
set of a leaf in a foliation of a cusped manifold. Obviously, the foliation itself
must be described in some computationally effective manner. Foliations carried by
Lackenby's taut ideal triangulations \cite{Lacktit} are good candidates for study.}
\end{enumerate}
\end{rmk}

\subsection{Realization questions}

\begin{qn}[Thurston]\label{quasigeodesic}
Suppose $M$ an atoroidal $3$--manifold admits an essential lamination. Does it
admit a (necessarily genuine) lamination with quasi--geodesic leaves?
\end{qn}

\begin{rmk}
\begin{enumerate}
\item{Any atoroidal manifold with an essential lamination
contains a genuine lamination $\Lambda$, and therefore has word--hyperbolic
fundamental group (\cite{Calp3},\cite{GK98}).
If such a lamination is not quasi--geodesic, perhaps one can find another
(more quasi--geodesic?) lamination transverse to it which is an ``eigenlamination''
for the extrinsic distortion of the geometry of $\til{\Lambda}$. As a caveat, one
should note that in general there is no naive notion of a {\em universal circle} for
a genuine lamination, so some new ideas are necessary. Note that faithful
circle actions are associated to certain essential laminations in \cite{CD01},
but these do not give monotone relations between circles at infinity of leaves,
as a real universal circle should. Rather they function as a kind of
circular parameterization of the {\em leaf space} of the laminations. Moreover, even
such naive circular parameterizations do not exist for certain non--tight genuine
laminations.}
\item{A lamination with quasi--geodesic leaves is tight, so a first step would be
to try to replace a non--tight lamination with a tight one. Associated to a
cataclysm (i.e. a collection $\lbrace \mu_i \rbrace$ of incomparable 
non--separated leaves  of $\til{\Lambda}$
which are the limit of a comparable monotone sequence) in a non--tight lamination
$\Lambda$, there is a geodesic lamination of $\H^2$ whose complementary regions
parameterize the $\mu_i$.}
\item{One reason to be interested in the existence of such a lamination is that it
would be an important technical advantage in
Thurston's program to geometrize laminar $3$--manifolds. Given
a choice of conformal structure on the each pair
$$(\text{complementary region } C, \text{boundary leaf of }C)$$
 leaves of $\Lambda$, one can try to generalize
the {\em skinning map} as a map from $\teich(\bigcup_C \partial C)$ to itself. In the case
of a taut foliation, it might proceed as follows. A universal
circle gives a monotone relation between the circles at infinity of the leaves
of $\til{\F}$. In fact, for $\lambda_i$ a maximal collection of incomparable leaves
above or below $\mu$, for any leaf $\mu$ of $\til{\F}$, there is a monotone relation
between $S^1_\infty(\mu)$ and a ``cactus'' made from the the union
$\bigcup_i S^1_\infty(\lambda_i)$. For each leaf $\lambda$ of $\til{\F}$,
we should choose a pair of points $$\tau^\pm(\lambda) \in \text{Teich}(D^2) \times
\text{Teich}(\overline{D^2})$$
a product of universal Teichm\"uller spaces as an initial
guess. Then monotone relations between circles at infinity give a new
marking $\sigma(\tau^+)$ obtained by sewing together the markings $\tau^+(\lambda_i)$
by the monotone relation. There are numerous problems with such a program;
for instance, there are many possibilities for leaves $\lambda_i$; moreover, there is
no natural decomposition of a laminated manifold along a collection of leaves into
a compact ``fundamental domain'', especially if $\F$ is minimal. Secondly, the
domain obtained by gluing together copies of $\tau^+(\lambda_i)$ will typically
not be a quasidisk, but will be pinched along an infinite lamination. One can
try to remedy this problem by letting $\sigma(\tau^+)$ be equal, not to the
union of the $\tau^+(\lambda_i)$ sewn together, but by $\tau^+$ bent in the direction
of this pinching lamination. But more importantly,
no foliation of a hyperbolic manifold can have all leaves quasi--geodesically
embedded; in particular, for all but the most trivial cases, there should be no
fixed point for the skinning map in 
$(\text{Teich}(D^2)\times \text{Teich}(\overline{D^2}))^L$, but rather in some
natural compactification, analogous to the Thurston boundary of Teichm\"uller space.

For a quasigeodesic lamination, at least a fixed point for the skinning map should
exist in some universal Teichm\"uller space, and proving geometrization in this
case might be considerably easier.}
\end{enumerate}
\end{rmk}

\subsection{Interactions of foliations with hyperbolic structures}

\begin{qn}\label{short_geodesics}
What do short geodesics look like with respect to taut foliations? Is there a
universal $\epsilon$ such that for every hyperbolic manifold $M$, every taut
foliation $\F$ of $M$, and every geodesic $\gamma$ with $|\gamma|<\epsilon$,
$\gamma$ is either isotopic into a leaf of $\F$ or isotopic to be transverse to $\F$?
What about homotopic?
\end{qn}
\begin{rmk}
\begin{enumerate}
\item{For a fixed taut foliation $\F$, leaves of $\til{\F}$ are uniformly
properly embedded in their $\epsilon$--neighborhoods, for some uniform
$\epsilon$. If $\F$ does not have $2$--sided branching, {\em every} loop is
homotopic to either a transverse or a tangential loop; this question asks whether
branching of $\til{\F}$ takes place on a scale uniformly bounded below over all
taut foliations.}
\item{For $M$ a hyperbolic manifold with a short geodesic $\gamma$, it is an
empirical observation of Rubinstein, Jaco and Sedgewick that a $1$--efficient
triangulation (see \cite{RuJa}) has a ``layered'' solid torus corresponding to
$\gamma$, which has a very standard combinatorial form. We can try to normalize
$\F$ relative to such a triangulation, possibly evacuating a subfoliation in
the process (\cite{Gab99}). Can one use the explicit combinatorial structure
of the layered solid torus to control the interaction of $\gamma$ with $\F$?}
\end{enumerate}
\end{rmk}

\begin{qn}
Is there a uniform bound on the Godbillon--Vey invariants of the taut foliations
of a hyperbolic manifold in terms of its volume?
\end{qn}

\rmk Some estimate for foliations on fixed $M$
can be obtained from the theorem of Gabai \cite{Gab99}
that there is a fixed triangulation $\tau$ of $M$ such that
every taut foliation of $M$ can be
made normal with respect to $\tau$. On the other hand, there are no known
{\it a priori} bounds on how many simplices such a $\tau$ must have in terms of the
number of simplices in a minimal triangulation,
or what relation $\tau$ has to the geometry of $M$. For instance, there
are geodesic triangulations which cannot be made ``virtually transverse'' to certain
taut foliations; i.e. there is no finite cover in which the pullback of the
triangulation can be normalized relative to the pullback of the foliation.

\begin{qn}
Suppose $\F$ is a taut foliation of a hyperbolic $3$--manifold $M$.
Let $$\pi:\til{M} \to L$$
be the projection to the leaf space of $\til{\F}$.
\begin{enumerate}
\item{For $\gamma$ a random walk in $\til{M}$
(which is isometric to $\H^3$), what is the typical behaviour of $\pi(\gamma)$?}
\item{Does a random walk in $\til{M}$ converge to a definite end of $L$?}
\item{What if we replace ``random walk'' with ``random geodesic'' in the previous
question?}
\item{Is the pushforward of asymptotic behaviour well--defined? That is, is it
true that for a set of geodesics $\gamma$ of full measure, for all $\gamma'$
a bounded distance from some $\gamma$ the behaviour of the pushforward of
$\gamma$ and $\gamma'$ exit the same end of $L$?}
\item{Suppose $\F$ has one--sided branching. Does a random walk always
exit $L$ in the unbranching direction?}
\item{Suppose there is a positive probability of a random walk exiting
a proper positive or negative end of $L$. Does this imply $\F$ is $\R$--covered?
What about if there is a positive probability of a random walk being recurrent?}
\end{enumerate}
\end{qn}

\begin{rmk}
\begin{enumerate}
\item{For the basic theory of random walks and Brownian motion, see e.g. \cite{Strprob}.}
\item{For $\F$ transversely measured, $L$ is naturally isomorphic to $\R$. Then the
pushforward of a random walk or geodesic
on $\til{M}$ under $\pi$ is a random walk on $\R$, and in particular the
pushforward is recurrent.}
\item{The question about foliations with one--sided branching has a positive
answer if the limit set of some (and therefore every) leaf of $\til{\F}$ is
not all of $S^2_\infty$. For in this case, the limit set of any leaf
is a dendrite of measure $0$, and the complement is entirely on the
nonbranching side of the leaf.}
\item{One has to be very careful in defining the behaviour of a random walk on
$L$; coarsely equivalent models for random walks on $\til{M}$ (e.g. combinatorial,
piecewise geodesic, etc.) might give rise to very different behaviour in $L$.
For instance, for many $\F$ with $2$--sided branching, we can find a pair of proper
lines $\gamma_1,\gamma_2$ in $\til{M}$ which are asymptotic (i.e. the distance
between them converges to $0$) but for which $\gamma_1$ exits an end of $L$ and
$\gamma_2$ is contained in a leaf. On the other hand, if $\gamma_1,\gamma_2$
are a bounded distance apart, and $\gamma_1$ exits to an {\em incomparable end}
of $L$ (i.e. the limit of a sequence of leaves $x_1 > y_1 < x_2 > y_2 < x_3 > y_3 < \dots$ where the $x_i$ are all incomparable and the $y_i$ are all incomparable) then
$\gamma_2$ exits to the {\em same} incomparable end. Moreover, if $\F$ arises from
a slithering or branched slithering, if $\gamma_1$ exits {\em any} end, $\gamma_2$
a bounded distance away exits the same end. If $\pi(\gamma)$ exits an end of $L$,
with probability one does it exit an incomparable end? (Yes in many cases)}
\item{If $\F$ is finite depth, we can ``trap'' the behaviour of a random walk
by compact leaves. If $\F$ is not a surface bundle over a circle, the random
walk has a definite probability of making an incomparable detour in a fixed
time interval. The tree of incomparable detours branches with ramification
bounded below; so a random walk will exit an incomparable end with probability $1$
in this case. Probably the same result is true if $\til{\F}$ contains any leaf which
is asymptotically separated.
Conjecturally, every leaf of every $\F$ with $2$--sided branching should
be asymptotically separated. However, part of the motivation for studying this
question would be to address this conjecture, rather than the other way around.}
\end{enumerate}
\end{rmk}

\begin{qn}
Suppose $\Lambda$ is an essential lamination of a hyperbolic manifold $M$.
Is $\Lambda$ isotopic to a lamination whose curvature is bounded below everywhere
by $-2$?
\end{qn}

\begin{rmk}
A complete properly embedded minimal surface in $\H^3$ has curvature bounded
below by $-2$ everywhere. Also see question~\ref{minimal_surface}.
\end{rmk}

\section{Foliated Teichm\"uller theory}

\subsection{Deformations, mapping--class groups}

For $\F$ a taut foliation of an atoroidal $3$--manifold, any conformal class of leafwise
metric on $\F$ can be ``uniformized'' to a leafwise hyperbolic metric. Co--deformations of such
metrics are parameterized by leafwise quadratic holomorphic differentials, and
presumably analogues of many other elements of one--dimensional complex analysis can be
found.

\begin{qn}
What kind of nontrivial ``mapping class elements'' are possible for
taut foliations?
\end{qn}

\begin{rmk}
\begin{enumerate}
\item{If $\F$ is a foliation of an atoroidal manifold $M$ arising
from a slithering over $S^1$, there is a canonical (up to
isotopy) transverse pseudo--Anosov flow $X$ associated to the slithering.
This flow is not in general leaf preserving at all times,
but the time one flow takes $\F$ to itself, and qualifies as a ``mapping class
element'' of infinite order. Slitherings, and slightly more generally, uniform
foliations are the only foliations for which
the time one flow of a transverse vector field takes leaves to leaves; are all
such mapping class elements isotopic to
the canonical automorphism?}
\item{For arbitrary taut $\F$ and atoroidal $M$, the mapping class group of $M$ is
finite, by results of 
\cite{Calp3} and \cite{GK98b}. If $\phi$ is a mapping class of $M,\F$
of infinite order, for some finite power $\phi$ will be isotopic to the identity.
Lift $\phi$ to $\til{\phi}:\til{M} \to \til{M}$.
It follows that there is a universal $t$ such that for every leaf $\lambda$ of
$\til{\F}$, the Hausdorff distance between $\lambda$ and $\til{\phi}(\lambda)$
is $\le t$. In particular, the space of leaves between $\lambda$ and
$\til{\phi}(\lambda)$ is isotopic to a product, and so if $\F$ branches,
every leaf corresponding to a non--Hausdorff point of $L$ is fixed by $\til{\phi}$.
For example, if $\F$ is finite depth, the only leaves which can move are those
of highest depth. If $\F$ is $\R$--covered but not uniform, there is
a monotone equivalent foliation $\G$ to which an element homotopic to
$\phi$ descends, which fixes $\til{\G}$ leafwise.}
\end{enumerate}
\end{rmk}

\subsection{Transverse Teichm\"uller flows}

\begin{qn}
A foliation is taut iff it admits a volume--preserving transverse flow.
Pseudo--Anosov flows are good candidates for ``best'' such transverse flows,
when they exist, which
is frequently.
Is there an analytic construction of pseudo--Anosov flows, by analogy with
Bers' proof of Thurston's classification of surface automorphisms?
\end{qn}

\rmk Fix a conformal sturcture on $\F$, and let $C$ be a transverse cone field.
Amongst unit speed vector fields $X$ transverse to $\F$, there are quasiconformally
extremal ones, by a standard compactness argument. It is not {\it a priori}
clear that there should be an extremal vector field amongst {\em all} transverse
unit speed vector fields. Moreover, one would want a ``best'' flow to
be metric independent. If $\F$ arises from a slithering, it makes sense to
parameterize a transverse vector field so that the time $1$ flow is a unit of the
slithering; the dilatation of the time $1$ map is well--defined with respect
to a leafwise conformal structure. Of course, one must simultaneously optimize over
all transverse flows {\em and all choices of leafwise conformal structure}.

For more general foliations, one needs to decide how to naturally parameterize
the leafwise flow in order to measure (infinitesimal) dilatation. One
plausible choice is to measure the result of flow with (transverse)
speed equal to the magnitude of a harmonic transverse measure.

\begin{qn}
Suppose $M$ is atoroidal and $\F$ arises from a slithering over $S^1$. Let $X$
be pseudo--Anosov transverse to $\F$, such that the time $1$ flow $Z$ takes
$\F$ to itself. Lift to $\til{M}$ and let $\lambda, Z^n(\lambda)$ be leaves of $\til{\F}$,
both uniformized as $\H^2$ by Candel's theorem.
Can $Z$ be approximated by mapping class elements between compact surfaces? That is,
are there integers $n_i$, a sequence $\Sigma_i$ of hyperbolic surfaces and
$\phi_i:\Sigma_i \to \Sigma_i$ Teichm\"uller representatives for the 
isotopy class $[\phi_i]$ such that the composition 
$$\til{\phi_i^{-1}}Z^{n_i}:\H^2 \to \H^2$$
is a $k_i$--quasi--isometry, where $k_i \to 1$?
\end{qn}
\begin{rmk}
\begin{enumerate}
\item{Implicitly in this question we are choosing an identification of $\H^2$
with $\lambda$ and $Z^n(\lambda)$. We could do this by choosing some $\alpha \in \pi_1(M)$
for which $\alpha Z^n(\lambda)$ is very close to $\lambda$ and compose it with a
nearest point map, or we could minimize $k_i$ over all isometries of $\H^2$ to itself.
Either method gives the same answer.}
\item{One reason to be interested in the existence of such $\phi_i$ is that it
might help to uniformize $M$.
We know that the mapping cylinders of the $\phi_i$ are hyperbolic, and therefore
estimates as above would give us approximate hyperbolic structures on larger and
larger subsets of $\til{M}$. What is lacking is an {\it a priori} comparison
between the geometry of a singular solv metric on the universal cover of a $3$--manifold
which fibers over $S^1$ and the geometry of the hyperbolic metric, depending only
on a bound on the degree of branching on the singular locus, and a
bound $c$ on the supremum of the distance from a point in $\til{M}$
to the singular locus. Note that without
such a bound, no {\it a priori} geometric comparison is possible.}
\end{enumerate}
\end{rmk}

\subsection{Algebraic and analytic geometry of foliations}

\begin{qn}
If $\F$ is a taut foliation, one can let $\gamma_i$ be a collection of transverse
circles to $\F$ intersecting every leaf and study the space of
functions $\O(\sum n_i\gamma_i)$ which are leafwise holomorphic, with poles
of order at most $n_i$ along $\gamma_i$. (Here the notation $\sum n_i\gamma_i$ is
meant to suggest a divisor on a Riemann surface). How do these function spaces
change as a function of $\gamma_i$? What is the effect of certain ``topological''
operations on the $\gamma_i$; e.g. crossing changes, cabling etc.
\end{qn}

\rmk Ghys used the fact that $\O(D)$ is large for certain $D$ to show that
taut foliations admit leafwise embeddings in complex projective space of sufficiently
high dimensions.

\section{Coarse foliations}

\begin{qn}
Suppose $\rho:\pi_1(M) \to \R$ is a $1$--cochain with bounded coboundary; i.e.
there is a uniform $C$ so that
$$|\rho(\alpha) + \rho(\beta) - \rho(\alpha\beta)|<C$$
for all $\alpha,\beta \in \pi_1(M)$. Consider $\pi_1(M)$ as a metric
space by thinking of it as the vertices of some Cayley graph.
Let $L_\rho = \rho^{-1}(I) \subset \pi_1(M)$.
\begin{enumerate}
\item{Are the coarse connected components of $L_\rho$ coarsely simply connected?}
\item{If the answer is ``yes'', is there a sense in which $L_\rho$ and its
translates by $\pi_1(M)$ can be thought of as coarse minimal planes in $\pi_1(M)$?}
\end{enumerate}
\end{qn}

\begin{rmk}
\begin{enumerate}
\item{If $L_\rho$ is coarsely connected for some $I$,
we say $\rho$ is {\em weakly uniform}, and
if it is both coarsely connected and coarsely simply connected, it is {\em uniform}.
One should think of $L_\rho$ as a coarse kernel for $\rho$; under this analogy,
coarsely connected corresponds to finitely generated, and coarsely simply connected
corresponds to finitely presented.

As one might guess from the Scott core theorem, for
$\pi_1(M)$ a $3$--manifold group, $\rho$ is weakly uniform iff it is uniform
(\cite{Calp4}). This may be thought of as a kind of ``coarse Stallings fibration
theorem'', or as a coarsening of Novikov's theorem. In \cite{Calp4} it is also
shown that such $\pi_1(M)$ either contain $\Z \oplus \Z$, or are word hyperbolic.}
\item{The potential usefulness of this question
comes from the fact that $1$--cochains with bounded coboundary
are very common on hyperbolic
$3$--manifolds. In fact, if the geometrization conjecture is true, {\em every}
$3$--manifold with infinite $\pi_1$ admits infinitely many independent such $\rho$.
This is in stark contrast to the fact that many $3$--manifolds with infinite
$\pi_1$ --- even hyperbolic $3$--manifolds --- do not admit taut foliations,
by \cite{RSS}.}
\end{enumerate}
\end{rmk}

\begin{qn}
Does every hyperbolic $3$--manifold admit a {\em taut cone field}? That is,
a cone field $C$ which is recurrent and supports only homotopically essential
loops.
\end{qn}
\begin{rmk}
\begin{enumerate}
\item{Taut cone fields are introduced in \cite{Calp2}. They exist on small Seifert
fibered spaces with infinite $\pi_1$ which do not admit essential laminations.
Examples of taut cone fields are cone fields supported by pseudo--Anosov flows.
However, not every hyperbolic $3$--manifold admits a pseudo--Anosov flow.
For instance, the Weeks manifold does not admit such a flow (\cite{CD01}).}
\item{Do cone fields stay taut after most surgeries
(i.e. those outside finitely many slopes) on a supported loop? Is there a class
of examples where the cone field stays taut for all nontrivial surgeries? See also
question~\ref{delman}.}
\item{By Sullivan \cite{Sul76}, there is a closed $2$--form positive on the de Rham
currents supported by a taut cone field. So, for instance, if $M$ is a rational
homology sphere, any complete surface (possibly immersed and noncompact)
transverse to the cone field enjoys a uniform linear isoperimetric inequality.}
\end{enumerate}
\end{rmk}

\subsection{Coarse invariants; deformation}

\begin{qn}
What deformations of a foliation or lamination should be thought of as
``inessential''? For instance --- monotone equivalence, cut--and--shear
along a surface or transverse lamination, isotopy of branch locus in
a branched cover, isomorphic universal circles etc.
\end{qn}

\rmk This question is prompted by the phenomenon that certain geometric
objects associated to foliations --- e.g. transverse pseudo--Anosov flows ---
are invariant under quite drastic deformations of the foliation.

An added bonus of a ``good'' deformation--equivalence class for foliations
might be a finite--dimensional parameter space, perhaps analogous to the
unit ball of the Thurston norm for finite--depth foliations. It seems that
an advantage of the theory of tight contact structures over the theory of
taut foliations is a better and more robust equivalence relation. E.g.
finiteness up to isotopy of tight contact structures on an atoroidal $3$--manifold
\cite{HG}.

\section{Numerical invariants}
\subsection{Godbillon--Vey}

The Godbillon--Vey invariant was introduced in \cite{GV71}. It is a invariance
of the foliated cobordism class of a sufficiently smooth ($C^2$ is sufficient)
codimension one foliation of a $3$--manifold.
A basic reference for this invariant is \cite{TsuGod}.

\begin{qn}
Suppose $\F$ is a minimal taut $C^2$ foliation of an atoroidal $3$--manifold $M$
with $$\gv(\F)[M] \ne 0$$ Is there a choice
of $1$--form $\alpha$ with $T\F = \ker(\alpha)$ for which the Godbillon--Vey
form $\omega$ (where $d\alpha = \alpha \wedge \omega$) has the same sign? i.e.
either $\omega \wedge d\omega \ge 0$ everywhere or $\le 0$ everywhere.
Say that such a foliation has {\em monotone wobble}.
\end{qn}

\begin{rmk} 
\begin{enumerate}
\item{Some condition on tori is essential. If $T$ is a torus leaf separating $M$
into two components $M_1 \coprod M_2$, then $\int_{M_i} \omega \wedge d\omega$
is independent
of the choice of $\alpha$, for each $i$. So we could splice together two
foliated manifolds with torus boundary, one with positive Godbillon--Vey
invariant, the other with negative invariant in such a way that the sum was
nonzero. More generally, if $M_1,M_2$ are two manifolds foliated by $\F_1,\F_2$,
and $\gamma_1,\gamma_2$ are transverse circles, the effect of {\em torus connect sum}
of $M_1, M_2$ along the $\gamma_i$ gives a foliated manifold whose Godbillon--Vey
invariant is the sum of the invariants for $M_1$ and $M_2$. A form $\alpha$ for
which $\omega$ is contact must have $\|\omega \wedge d\omega\|$ potentially very
large in a neighborhood of the separating torus
$\partial N(\gamma_1) = \partial N(\gamma_2)$.}
\item{Minimality is also essential, or else regions of positive and
negative twisting can be separated by some exceptional minimal set or even
closed leaf. On the other hand, for a minimal foliation, positive twisting
can be propagated along neighborhoods of a leafwise path to regions of
negative twisting, by analogy with the process of deforming a transitive confoliation
to a contact structure.}
\item{If $\omega \wedge d\omega \ge 0$ can be achieved everywhere, when can
it be achieved so that the zero locus $\omega \wedge d\omega = 0$ is a collection
of knots, and $M,\omega$ is locally modelled on a branched cover of a contact
manifold over a Legendrian link?}
\item{If $\omega \wedge d\omega >0$ everywhere --- that is, if $\omega$ is actually
a contact structure, then in particular $T\F$ admits a nowhere zero section, and
the Euler class of $T\F$ is zero. This implies that $\F$ has no closed leaves.}
\end{enumerate}
\end{rmk}

\begin{qn}
For $\F$ as in the previous question,
suppose there is a choice of $\alpha$ for which $\omega$ is a contact form.
Is the contact structure defined by $\omega$ necessarily tight?
\end{qn}

\subsection{Foliated Gromov norms}

For $M,\F$ a foliated manifold, the {\em foliated Gromov norm} on $H_i$ is the
infimum of the $L_1$ norm on chains representing homology classes, where we restrict
attention to chains which are {\em transverse} to $\F$ (in the sense that the
induced foliation of each simplex should be affine). A reference for this
material is \cite{Cal01}. Notice that to define the Godbillon--Vey invariant of
$\F$, the transverse structure must be at least Lipschitz $+$ quadratic variation 
(see \cite{Tsu95}); by contrast, the foliated Gromov norm depends only on
the {\em homeomorphism} type of the foliation.

\begin{qn}
Calculate the norm of the fundamental class of a hyperbolic $3$--manifold for some
taut foliation $\F$ with two--sided branching.
\end{qn}

\rmk If $\F$ is $\R$--covered or has one--sided branching, the foliated Gromov norm agrees
with the usual Gromov norm.

If some leaf of $\til{\F}$ is {\em asymptotically separated} 
(see question~\ref{asymptotically_separated}) then the value of the norm on this fundamental
class is {\em strictly} greater than the value of the usual Gromov norm on this class.

\begin{qn}
Let $\F,\G$ be taut foliations on a hyperbolic manifold $M$. Are there examples
where there is a finite cover of $M$ such that a sequence of isotopies of the lift
of $\G$ converges geometrically to the lift of $\F$, but no such sequence of
isotopies exists in $M$?  
\end{qn}

\rmk If this can be done, the foliated Gromov norm of $\F$ is at least as
large as that of $\G$ for every homology class.

\begin{qn}
Suppose $\F$ is a foliation (possibly $\R$--covered) of a hyperbolic $3$--manifold.
Define a foliated Gromov norm using {\em cubical} chains. Is the value of the foliated norm on
the fundamental class always {\em strictly} greater than the value of the usual (cubical)
Gromov norm?
\end{qn}

\rmk Cubical chains cannot necessarily be straightened rel. vertices to be transverse to
a given foliation. Suppose $\F$ is a transversely measured foliation. Lift this
transverse measure to a function $f:\til{\F} \to \R$. Then for a random regular cube $C$
of side length $t$, the ordering of the vertices of $C$ induced by $f$ converges to a
random distribution as $t \to \infty$.

\begin{qn}
What kinds of local order structure are there on a family of deformations
of a (taut) foliation? Can one use such structures to define co--ordinates
on the ``space of deformations'' of a taut foliation?
\end{qn}
\begin{rmk}
\begin{enumerate}
\item{For instance, local order structure could be provided by numerical
invariants, such as the Godbillon--Vey invariant, or foliated Gromov norms.}
\item{A family $\F_t$ obtained by cut and shear along some particular
submanifold, where the re--gluings $f_t$ vary monotonically in some sense
should be considered a monotone deformation. E.g. in local co--ordinates,
the differences $f_t - f_s$ should be non--negative functions for $t>s$.}
\end{enumerate}
\end{rmk}

\subsection{Invariants for laminations}

\begin{qn}
Is there some notion of a Godbillon--Vey invariant for a lamination?
\end{qn}

\rmk One can look
at eigendirections of pinching for holonomy and construct a leafwise vector field;
a Godbillon--Vey type invariant should be a topological measure of how much this
vector field twists in the positive direction. Even an invariant taking values in
$+,-,0$ would be interesting, where the invariant would be $+$ if for some metric,
the transverse twisting (``helical wobble'' is Thurston's term) is always positive,
$-$ if for some metric it is always negative, and $0$ if neither can be achieved.
Is this even well--defined? Maybe one would need to restrict attention to
the behavior at infinity --- e.g. ask for a metric such that for
some choice of partition into guts and interstices, the helical wobble over the
interstices is always positive.

\section{Immersed objects}

\begin{qn}
Is there a geometric notion for a $3$--manifold analogous to LERFness for foliations?
What properties could a manifold have so that immersed essential laminations are virtually embedded?
\end{qn}

\rmk It might be better to restrict attention to foliations dual to one--cochains
with bounded coboundary, since these have a more ``algebraic'' description.

\begin{qn}
Let $\F$ be a taut foliation of $M$. Can leaves of
$\F$ be approximated by compact essential surfaces? That is, given a leaf $\lambda$ of
$\F$ and a point $p \in \lambda$, is there a sequence of immersed incompressible
surfaces $\phi_i:\Sigma_i \to M$ and points $p_i \in \Sigma_i$ such that the
images under $\phi_i$ of the balls of radius $r_i$ about $p_i$ where $r_i \to \infty$
converge on compact sets to $p,\lambda$?
\end{qn}

\begin{rmk}
\begin{enumerate}
\item{If one does not ask for incompressibility, it is easy to find such
approximating surfaces; for, if $D \subset \lambda$ is a big embedded subsurface,
we can just take two parallel copies of $D$ tubed together by a collection of thin annuli,
one for each boundary component of $D$.}
\item{If $\pi_1(M)$ is LERF, each of the surfaces $\Sigma_i$ lifts to an embedded surface
in some finite cover. An interesting additional constraint is to ask for these lifts to
be Thurston norm minimizing.}
\item{Suppose we can choose the $\Sigma_i$ which converge {\em uniformly} to some leaf 
$\lambda$ in the sense that for every $\epsilon$ there is an $i$ such that 
$(\phi_i)_*T\Sigma_i$ makes an angle of less than $\epsilon$ at
every point with $T\lambda$. Then the de Rham cycles $\phi_i([\Sigma_i])$
converge to a de Rham cycle supported by $\F$; since $\F$ is taut,
this implies $H_1(M;\R)$ is nontrivial. On the other hand, the
support of the limiting cycle need not be all of $\F$. This can happen,
for instance, if the $\Sigma_i$ are surfaces of a fibration of $M$ over $S^1$ representing
homology classes projectively limiting to a class contained in a high codimension face
of the Thurston norm. However, if $\lambda$ is dense, every transverse loop to $\F$
is eventually transverse to some $\Sigma_i$, and therefore pairs
positively with the projective homology class represented by $\Sigma_i$.
It follows that every transverse loop to $\F$ is nontrivial in homology, and therefore
$\F$ is a perturbation of a surface bundle over $S^1$.}
\end{enumerate}
\end{rmk}

\subsection{Total foliations; families of foliations}

A {\em total} foliation on an $n$ manifold is a collection of $n$ foliations $\F_i$ for
$i = 1 \dots n$ which
are mutually transverse everywhere and in sufficiently small charts, are
topologically conjugate to the total foliation of $\R^n$ by translates of the
co--ordinate planes.

\begin{qn}
What $3$--manifolds admit total taut foliations?
\end{qn}

\rmk Every $3$--manifold admits a total foliation, by a local construction
due to Hartdorp \cite{Har}. Suspensions
of Anosov automorphisms of tori admit total taut foliations. Some examples
arise in the theory of slitherings (\cite{Thu97}). Anosov flows give rise to some
good examples. It is
probably better in general for $M$ a hyperbolic manifold to look for taut foliations of $M$
with {\em genuine} transverse pseudo--Anosov flows; these exist in very many examples.

\begin{qn}
Are there any interesting examples of total genuine laminations?
\end{qn}

\begin{rmk}
\begin{enumerate}
\item{One can generalize total foliations to a $\G$--structure where $\G$ is
the pseudogroup of homeomorphisms between subsets of $\R^3$ which permute translates
of the co--ordinate planes, but are not required to take horizontal to
horizontal, etc. For $M$ a manifold with a $\G$--structure, there is a
representation $\pi_1(M) \to {\sf S}_4$ whose kernel defines a finite cover
with a total foliation (in the usual sense). For {\em total laminations}, there
are local obstructions, much as there are local obstructions to co--orientability.
It is reasonable to use this more flexible definition of total lamination, since
it encompasses many known examples (e.g. pseudo--Anosov flows transverse to
foliations). Is it useful to think of a total genuine lamination as some
kind of intrinsically three--dimensional analogue of a pseudo--Anosov flow?}
\item{The immersed meridional surface in a cubulation of strictly negative curvature is
an example of a total genuine lamination. The complementary regions to a total
genuine lamination should be cubes, products of negatively curved polygons with intervals,
and polyhedra satisfying the condition of Andre'ev's theorem (see \cite{HodRiv}) for
all right--angles. It should be straightforward to find a direct proof that manifolds
with total genuine laminations admit CAT $0$ metrics.}
\item{In \cite{Li4plane} examples are constructed of hyperbolic manifolds which
do not admit nonpositive cubulations. It would be interesting to know whether
these manifolds can admit total genuine laminations.}
\end{enumerate}
\end{rmk}

\subsection{Amenability}

\begin{qn}
What is the weakest useful $2$--dimensional object that might be present in
every atoroidal $3$--manifold? For instance,
does every hyperbolic $3$--manifold $M$ contain an immersed quasigeodesic
surface $\Sigma$ of amenable growth?
\end{qn}

\begin{rmk} 
\begin{enumerate}
\item{Of course, the virtual infinite positive Betti number
conjecture would imply that every hyperbolic
$3$--manifold contains an immersed quasigeodesic {\em compact} surface. The point here
is to think of as weak a condition as possible, so that there are (hopefully?)
many examples to study.}
\item{The closure of the image of such an immersed surface $\Sigma$ supports a 
de Rham cycle by the usual method of Goodman--Plante \cite{GP79}. If $M$ is a rational
homology sphere, then for any loop $\gamma$ in $M$, large subsets of $\Sigma$
must intersect $\gamma$ positively and negatively with approximately equal frequency.
Does this imply ergodicity of the geodesic flow?}
\end{enumerate}
\end{rmk}

\section{Miscellaneous}

\begin{qn}\label{well_ordered}
What possibilities are there for (co--oriented) laminations in a $3$--manifold whose
transverse spaces are {\em well--ordered}?

Is there a (useful) theory of branched surfaces with ordinal--valued weights?
\end{qn}

\rmk
Probably the condition that laminations be geometrically realized in a $3$--manifold
implies that the order types will be {\em countable}. There is an
{\em evacuation} procedure to replace a well--ordered lamination by a simpler one:
given $\Lambda$, we can replace $\Lambda$ by $\Lambda_\infty$, the sublamination
whose leaves are labelled by limit ordinals in a given chart. This will be a proper
sublamination. To get an interesting theory, one would want to study laminations
for which this procedure does not strictly decrease some notion of complexity --- so
that the order types should be at least uncountable.

For branched surfaces with no realization issues, more exotic and interesting
ordinals could be used; perhaps one can use this to give ``geometric'' interpretations
of interesting theorems in logic.

\begin{qn}
Is there a good notion of {\em taut foliated cobordism}? Are there numerical invariants
of the equivalence classes this induces on taut foliations which are finer than the
Godbillon--Vey invariant?
\end{qn}

\end{document}